\documentclass[a4paper,12pt]{article}

 \usepackage{amsfonts,amsmath,mathrsfs,amssymb}
 \usepackage{times,helvet,courier,type1cm}

 \usepackage{harvard}
 \usepackage{color}

 \allowdisplaybreaks

 \newtheorem{thm0}{Theorem}[section]
 \newtheorem{exa0}{Theorem}[section]

 \newtheorem{def1}[thm0]{Definition}
 \newtheorem{lem1}[thm0]{Lemma}
 \newtheorem{thm1}[thm0]{Theorem}
 \newtheorem{cor1}[thm0]{Corollary}
 \newtheorem{pro1}[thm0]{Proposition}
 \newtheorem{con1}[thm0]{Condition}
 \newtheorem{rem1}[thm0]{Remark}
 \newtheorem{exa1}[exa0]{\it{Example}}

 \def\bgtheorem{\begin{thm1}}\def\edtheorem{\end{thm1}}
 \def\bgcorollary{\begin{cor1}}\def\edcorollary{\end{cor1}}
 \def\bgproposition{\begin{pro1}}\def\edproposition{\end{pro1}}
 \def\bgcondition{\begin{con1}}\def\edcondition{\end{con1}}
 \def\bgremark{\begin{rem1}}\def\edremark{\end{rem1}}
 \def\bgexample{\begin{exa1}\rm{}\def\edexample{\end{exa1}}}

 \def\benumerate{\begin{enumerate}}\def\eenumerate{\end{enumerate}}
 \def\bitemize{\begin{itemize}}\def\eitemize{\end{itemize}}\def\itm{\item}

 \def\beqlb{\begin{eqnarray}}\def\eeqlb{\end{eqnarray}}
 \def\beqnn{\begin{eqnarray*}}\def\eeqnn{\end{eqnarray*}}

 \def\eqref#1{{\rm(\ref{#1})}}

 \def\ar{\!\!\!&}\def\nnm{\nonumber}\def\ccr{\nnm\\}

 \def\<{\langle}\def\>{\rangle}

 \def\mcr{\mathscr}\def\mbb{\mathbb}\def\mbf{\mathbf}
 \def\mrm{\mathrm}

 \def\proof{\noindent{\textit{Proof.~~}}}
 \def\qed{\hfill$\square$\smallskip}

 \def\d{\mrm{d}}\def\e{\mrm{e}}
 \def\var{\mrm{var}}

 \def\supp{\mrm{supp}}

\begin{document}

\noindent{(Version: 2019-07-23)}

\bigskip\bigskip\bigskip

\centerline{\LARGE Ergodicities and exponential ergodicities of}

\centerline{\LARGE Dawson--Watanabe type processes\,\footnote{Supported by the
National Natural Science Foundation of China (No.11531001).}}

\bigskip

\centerline{\Large Zenghu Li\,\footnote{\textit{E-mail address:} lizh@bnu.edu.cn.}}

\medskip

\centerline{\it School of Mathematical Sciences, Beijing Normal University,}

\centerline{\it Beijing 100875, P.R.~China. }

\bigskip

{\narrower

\noindent{\textit{Abstract:}} Under natural assumptions, we prove the ergodicities and exponential ergodicities in Wasserstein and total variation distances of Dawson--Watanabe superprocesses without or with immigration. The strong Feller property in the total variation distance is derived as a by-product. The key of the approach is a set of estimates for the variations of the transition probabilities. The estimates in Wasserstein distance are derived from an upper bound of the kernels induced by the first moment of the superprocess. Those in total variation distance are based on a comparison of the cumulant semigroup of the superprocess with that of a continuous-state branching process. The results improve and extend considerably those of Stannat (2003a, 2003b) and Friesen (2019+). We also show a connection between the ergodicities of the associated immigration superprocesses and decomposable distributions.

\bigskip

\noindent{\textit{Key words and phrases.}} Dawson--Watanabe superprocess, immigration, coupling, strong Feller property, stationary distribution, exponential ergodicity.

\smallskip

\noindent{\textit{MSC 2010 Subject Classification.}} Primary 60J68; secondary 60J80, 60J35

\par}

\bigskip

\section{Introduction}

 \setcounter{equation}{0}

Measure-valued branching processes (MB-processes) have been studied extensively in the past decades. They arise naturally in the study of rescaling limits of branching particle systems. A special class of those processes are known as Dawson--Watanabe superprocesses. Measure-valued branching processes with immigration (MBI-processes) are generalizations of the MB-processes that consider the input into the system from outside sources. The reader may refer to Dawson (1993), Dynkin (1994), Etheridge (2000), Le~Gall (1999), Li (2011) and the references therein for the literature in the subject. When the underlying space is a finite set, the MB-process is often referred to as a continuous-state branching process (CB-process); see, e.g., Rhyzhov and Skorokhod (1970) and Watanabe (1969).

The strong Feller and ergodic properties are both important topics in the theory of Markov processes. In particular, a necessary and sufficient condition for the ergodicity of a one-type continuous-state branching processes with immigration (CBI-process) was announced in Pinsky (1972). A proof of the result can be found in Li (2011). It was proved in Keller-Ressel and Mijatovi\'c (2012) that the class of stationary distributions of one-type CBI-processes is strictly contained in the class of infinitely divisible distributions on the positive half-line and is strictly larger than that of classical self-decomposable distributions. The strong Feller property and exponential ergodicity of the process in the total variation distance were proved in Li and Ma (2015) using a coupling method; see also Li (2019+). The exponential ergodicity played an important role in the study of asymptotics of the estimators for the process in Li and Ma (2015).

A CBI-process involves affine structures in its generator and the Laplace transform of its transition probabilities in a similar way as an Ornstein--Uhlenbeck type (OU-type) process. The class of affine Markov processes unifies the CBI- and OU-type processes. This unified treatment of those processes has developed interesting connections between several areas in theory and applications of probability. A sufficient condition for the ergodicity of the OU-type process in the sense of weak convergence was given in Sato and Yamazato (1984). The coupling and strong Feller properties of those processes were studied in Wang (2011a) and some gradient estimates were given in Wang (2011b). Using coupling techniques, Schilling and Wang (2012) and Wang (2012) investigated the ergodicity and exponential ergodicity of the processes in the total variation distance. The strong Feller properties and exponential ergodicity of OU-type processes in Banach spaces were studied in Wang and Wang (2013). A result on exponential ergodicity of affine processes in Wasserstein distance was proved in the very recent work of Friesen et al. (2019+); see also Jin et al.\ (2018+).

The general immigration structures associated with a measure-valued branching process can be formulated in terms of skew-convolution semigroups. It was proved in Li (1996a, 2011) that such a semigroup is uniquely determined by an infinitely divisible probability entrance law. When the entrance law is closable, the immigration is governed by an infinitely divisible distribution. The study of exponential ergodicities of the Dawson--Watanabe superprocess with immigration in Wasserstein and total variation distances was initiated by Stannat (2003a, 2003b), who considered a Feller underlying process and a local branching mechanism. Stannat (2003a, 2003b) also focused on a particular immigration structure determined by a finite measure on the underlying space. The results on exponential ergodicity in Wasserstein distance of Stannat (2003a, 2003b) were generalized in the recent work of Friesen (2019+) to a Borel right underlying processes and a nonlocal branching mechanism.

The main purpose of this work is to study the ergodicities and exponential ergodicities of Dawson--Watanabe superprocesses without or with immigration in general settings. We shall prove those properties in Wasserstein and total variation distances under natural assumptions. The strong Feller property will be derived as a by-product. The key of the approach here is a set of estimates for the distances between the relevant distributions. The estimates in Wasserstein distance are derived from an upper bound of the kernels induced by the first moment of a superprocess. Those in total variation distance are based on a comparison of the cumulant semigroup of the superprocess with that of a one-dimensional CB-process. The approach is simpler than that of analysis of generators used in Stannat (2003a, 2003b). To illustrate the essential structures, we shall establish results for general MB- and MBI-processes and then specify them to the case of superprocesses without or with immigration. The immigration structures considered here are determined by infinitely divisible probability entrance laws not necessarily closable. The results improve and extend considerably those of Stannat (2003a, 2003b) and Friesen (2019+). In fact, we give accurate evaluations of the distances between some of the distributions. We also show that the ergodicities are closely related with some self-decomposable distributions.

The paper is organized as follows. In Section~2, the results for general MB-processes are presented. The results for Dawson--Watanabe superprocesses are given in Section~3. In Section~4, we study the immigration structures and related ergodicities. The connection between the ergodicities and self-decomposable distributions is explained in Section~5. In Section~6, we give some examples including comparisons of the results with those of Stannat (2003a, 2003b) and Friesen (2019+).

\section{General MB-processes}

 \setcounter{equation}{0}

Consider a Lusin topological space $E$, i.e., a homeomorph of a Borel subset of some compact metric space. Let $M(E)$ be the space of finite Borel measures on $E$ furnished with the topology of weak convergence. Then $M(E)$ is also a Lusin topological space; see, e.g., Theorem~1.16 of Li (2011). Let $B(E)$ be the Banach space of bounded Borel functions on $E$ equipped with the supremum norm $\|\cdot\|$. Let $B(E)^+\subset B(E)$ denote the subset of positive ($=$\,nonnegative) functions. For $\mu\in M(E)$ and $f\in B(E)$ write $\mu(f)=\int_Ef\d\mu$. A bounded kernel $\gamma$ on $E$ induces an operator on $B(E)$ defined by
 \beqnn
\gamma f(x)= \gamma(x,f)= \int_Ef(y)\gamma(x,\d y), \qquad x\in E, f\in B(E).
 \eeqnn
The kernel also induces an operator on $M(E)$ defined by
 \beqnn
\mu\gamma(f)= \int_E\mu(\d x)\int_Ef(y)\gamma(x,\d y), \qquad \mu\in M(E), f\in B(E).
 \eeqnn
It is well-known that a probability measure $Q$ on $M(E)$ is uniquely determined by its Laplace functional $L_Q$ defined by
 \beqnn
L_Q(f)= \int_{M(E)} \e^{-\nu(f)} Q(\d\nu), \qquad f\in B(E)^+.
 \eeqnn
For $\mu$ and $\nu\in M(E)$ let $|\mu-\nu|$ denote the \textit{total variation} of the signed-measure $\mu-\nu$. Then $\|\mu-\nu\|_{\var}:= |\mu-\nu|(E)$ is the \textit{total variation distance} between  $\mu$ and $\nu$. For a function $F$ on $M(E)$, its Lipschitz constant relative to the total variation distance  is defined by
 \beqnn
L_{\var}(F) = \sup\big\{\|\mu-\nu\|_{\var}^{-1}|F(\mu)-F(\nu)|: \mu\neq \nu\in M(E)\big\}.
 \eeqnn
A \textit{coupling} of two probability measures $Q_1$ and $Q_2$ on $M(E)$ is a probability measure $P$ on $M(E)^2$ with marginals $P(\cdot\times E)= Q_1(\cdot)$ and $P(E\times \cdot)= Q_2(\cdot)$. The \textit{Wasserstein distance} $W_1(Q_1,Q_2)$ between $Q_1$ and $Q_2$ is defined by
 \beqnn
W_1(Q_1,Q_2)= \inf_P\int_{M(E)^2} \|\mu-\nu\|_{\var} P(\d\mu,\d\nu),
 \eeqnn
where $P$ runs over all couplings of $Q_1$ and $Q_2$. We refer to Chen (2004a, 2004b) for systematic discussions of couplings and Wasserstein distances.

A conservative Markov process $X$ with state space $M(E)$ is called a \textit{measure-valued branching process} (MB-process) if its transition semigroup $(Q_t)_{t\ge 0}$ satisfies the (regular) \textit{branching property}:
 \beqlb\label{2.1}
\int_{M(E)} \e^{-\nu(f)} Q_t(\mu,\d\nu)
 =
\exp\{-\mu(V_tf)\}, \qquad \mu\in M(E), f\in B(E)^+,
 \eeqlb
where
 \beqnn
V_tf(x)= -\log\int_{M(E)} \e^{-\nu(f)} Q_t(\delta_x,\d\nu), \qquad x\in E.
 \eeqnn
By \eqref{2.1} we have $\mu(V_tf)< \infty$ for every $\mu\in M(E)$, so $V_tf\in B(E)^+$. It is not hard to show that the operators $(V_t)_{t\ge 0}$ on $B(E)^+$ satisfy $V_sV_t= V_{s+t}$ for $s\ge 0$ and $t\ge 0$. We call $(V_t)_{t\ge 0}$ the \textit{cumulant semigroup} of $X$. From \eqref{2.1} we see that $(Q_t)_{t\ge 0}$ has the \textit{branching property}
 \beqlb\label{2.2}
Q_t(\mu_1+\mu_2,\cdot)
 =
Q_t(\mu_1,\cdot)*Q_t(\mu_2,\cdot), \qquad t\ge 0, ~\mu_1\mu_2\in M(E),
 \eeqlb
where ``$*$'' denote convolution. By Theorem~2.4 of Li (2011) we have the L\'{e}vy--Khintchine type representation:
 \beqlb\label{2.3}
V_tf(x)= \lambda_t(x,f) + \int_{M(E)^\circ} \big(1-\e^{-\nu(f)}\big) L_t(x,\d\nu), \quad x\in E, f\in B(E)^+,
 \eeqlb
where $\lambda_t(x,\d y)$ is a bounded kernel on $E$ and $[1\land \nu(1)]L_t(x,\d\nu)$ is a bounded kernel from $E$ to $M(E)^\circ$.

From \eqref{2.1} we see that $0\in M(E)$ is a trap for $(Q_t)_{t\ge 0}$. Then the Dirac measure $\delta_0$ is a stationary distribution for $(Q_t)_{t\ge 0}$. Moreover, we have $\lim_{t\to \infty}Q_t(\mu,\cdot) = \delta_0$ by weak convergence for every $\mu\in M(E)$ if and only if $\lim_{t\to \infty}V_t1(x)\to 0$ for every $x\in E$.

\bgcondition\label{t2.1} For $t\ge 0$ and $f\in B(E)$, the following function is bounded on $E$:
 \beqnn
\pi_tf(x):= \int_{M(E)} \nu(f) Q_t(\delta_x,\d\nu), \qquad x\in E.
 \eeqnn
\edcondition

Under Condition~\ref{t2.1}, the MB-process with deterministic initial state has finite moments. In fact, by the branching property \eqref{2.2} it is not hard to show that the family of kernels $(\pi_t)_{t\ge 0}$ on $E$ constitute a semigroup and
 \beqlb\label{2.4}
\int_{M(E)} \nu(f) Q_t(\mu,\d\nu) = \mu(\pi_tf), \qquad \mu\in M(E), f\in
B(E).
 \eeqlb
By Jensen's inequality, we have $V_tf(x)\le \pi_tf(x)$ for $x\in E$ and $f\in B(E)^+$.

The next theorem gives upper and lower bounds for the variations in Wasserstein distance of the transition probabilities of the MB-process started from two different initial states.

\bgtheorem\label{t2.2} Suppose that Condition~\ref{t2.1} holds. Then for $t\ge 0$ and $\mu, \nu\in M(E)$ we have
 \beqlb\label{2.5}
|\mu(\pi_t1)-\nu(\pi_t1)|\le W_1(Q_t(\mu,\cdot),Q_t(\nu,\cdot))\le |\mu-\nu|(\pi_t1).
 \eeqlb
\edtheorem

\proof Let $F_1(\eta)= \eta(1)$ for $\eta\in M(E)$. Then $F_1$ is a Lipschitz function on $M(E)$ in the total variation distance with $L_{\var}(F_1)= 1$. By Theorem~5.10 in Chen (2004a, p.181) we have
 \beqnn
W_1(Q_t(\mu,\cdot),Q_t(\nu,\cdot))
 \ge
\int_{M(E)} \eta(1)(Q_t(\mu,\d\eta)-Q_t(\nu,\d\eta))
 =
\mu(\pi_t1)-\nu(\pi_t1).
 \eeqnn
Similarly we have $W_1(Q_t(\mu,\cdot),Q_t(\nu,\cdot))\ge \nu(\pi_t1)-\mu(\pi_t1)$. Then the first inequality in  \eqref{2.5} follows. Let $(\mu-\nu)_+$ and $(\mu-\nu)_-$ denote the upper and lower variations of the signed measure $\mu-\nu$ in its Jordan-Hahn decomposition, respectively. Let $\mu\land\nu= \mu - (\mu-\nu)_+= \nu - (\mu-\nu)_-$. Let $P_t(\mu,\nu,\d\gamma_1,\d\gamma_2)$ be the image of the product measure
 \beqnn
Q_t(\mu\land\nu,\d\eta_0) Q_t((\mu-\nu)_+,\d\eta_1) Q_t((\mu-\nu)_-,\d\eta_2)
 \eeqnn
under the mapping $(\eta_0,\eta_1,\eta_2)\mapsto (\gamma_1,\gamma_2):= (\eta_0+\eta_1, \eta_0+\eta_2)$. Then $P_t(\mu,\nu,\d\gamma_1,\d\gamma_2)$ is a coupling of $Q_t(\mu,\d\gamma_1)$ and $Q_t(\nu,\d\gamma_2)$. It follows that
 \beqnn
W_1(Q_t(\mu,\cdot),Q_t(\nu,\cdot))
 \ar\le\ar
\int_{M(E)^2}\|\gamma_1-\gamma_2\|_{\var}P_t(\mu,\nu,\d\gamma_1,\d\gamma_2) \cr
 \ar=\ar
\int_{M(E)}Q_t(\mu\land\nu,\d\eta_0) \int_{M(E)}Q_t((\mu-\nu)_+,\d\eta_1) \cr
 \ar\ar\qquad
\int_{M(E)}\|\eta_1-\eta_2\|_{\var}Q_t((\mu-\nu)_-,\d\eta_2) \cr
 \ar\le\ar
\int_{M(E)}Q_t(\mu\land\nu,\d\eta_0) \int_{M(E)}Q_t((\mu-\nu)_+,\d\eta_1) \cr
 \ar\ar\qquad
\int_{M(E)}[\eta_1(1)+\eta_2(1)]Q_t((\mu-\nu)_-,\d\eta_2) \cr
 \ar=\ar
\int_{M(E)}\eta(1)Q_t(|\mu-\nu|,\d\eta)
 =
|\mu-\nu|(\pi_t1),
 \eeqnn
where we have used the relation $|\mu-\nu|= (\mu-\nu)_+ + (\mu-\nu)_-$ and the branching property \eqref{2.2}. Then \eqref{2.5} follows. \qed

\bgcorollary\label{t2.3} Suppose that Condition~\ref{t2.1} holds. Then $L_{\var}(Q_tF)\le \|\pi_t1\|L_{\var}(F)$ for any $t\ge 0$ and Borel function $F$ on $M(E)$. \edcorollary

\bgcorollary\label{t2.4} Suppose that Condition~\ref{t2.1} holds. Then $W_1(Q_t(\mu,\cdot),Q_t(\nu,\cdot)) = (\mu-\nu)(\pi_t1)$ for $t\ge 0$ and $\mu\ge \nu\in M(E)$.
\edcorollary

\bgcorollary\label{t2.5} Suppose that Condition~\ref{t2.1} holds. Then $W_1(Q_t(\mu,\cdot),\delta_0)= \mu(\pi_t1)\to 0$ as $t\to \infty$ for every $\mu\in M(E)$ if and only if $\lim_{t\to \infty} \pi_t1(x)= 0$ for every $x\in E$.
\edcorollary

By Corollary~\ref{t2.5}, the class of Lipschitz functions on $M(E)$ in the total variation distance is invariant under the transition semigroup of the MB-process. To give some estimates for the variations in total variation distance of the transition probabilities of the process, let us consider the following condition:

\bgcondition\label{t2.6} For each $t>0$ the function $\bar{V}_t(x):= \lim_{\lambda\to \infty} V_t\lambda(x)$ is bounded on $E$.
\edcondition

\bgproposition\label{t2.7} Suppose that Condition~\ref{t2.6} holds. Then we have \eqref{2.3} with $\lambda_t(x,1)= 0$ and $\bar{V}_t(x)= L_t(x,M(E)^\circ)$ for $t>0$ and $x\in E$. Moreover, the mapping $t\mapsto \bar{V}_t(x)$ on is decreasing $(0,\infty)$ and
 \beqlb\label{2.6}
Q_t(\mu,\{0\})= \e^{-\mu(\bar{V}_t)}, \qquad t>0,\mu\in M(E).
 \eeqlb
\edproposition

\proof The first assertion is immediate. For $t\ge r>0$ and $x\in E$, by taking $\mu=\delta_x$ in \eqref{2.1} and using monotone convergence we have $\bar{V}_t(x)= \lim_{\lambda\to \infty}V_t\lambda(x)= \lim_{\lambda\to \infty}V_r V_{t-r} \lambda(x)= V_r\bar{V}_{t-r}(x)\le \bar{V}_r(x)$. From \eqref{2.1} we get \eqref{2.6}. \qed

The reader may refer to Dawson (1993, p.195) for an earlier form of the above result. The next theorem gives some general estimates for the variations in total variation distance of the transition probabilities of the MB-process.

\bgtheorem\label{t2.8} Suppose that Condition~\ref{t2.6} holds. Then, for $t>0$ and $\mu,\nu\in M(E)$,
 \beqlb\label{2.7}
2|\e^{-\mu(\bar{V}_t)}-\e^{-\nu(\bar{V}_t)}|
 \le
\|Q_t(\mu,\cdot)-Q_t(\nu,\cdot)\|_{\var}
 \le
2(1-\e^{-|\mu-\nu|(\bar{V}_t)}).
 \eeqlb
\edtheorem

\proof If $\mu(\bar{V}_t)\le \nu(\bar{V}_t)$, by \eqref{2.6} we have $Q_t(\mu,\{0\})-Q_t(\nu,\{0\})= \e^{-\mu(\bar{V}_t)} - \e^{-\nu(\bar{V}_t)}\ge 0$, and so $\|Q_t(\mu,\cdot)-Q_t(\nu,\cdot)\|_{\var}\ge 2(\e^{-\mu(\bar{V}_t)} - \e^{-\nu(\bar{V}_t)})$. Then the first inequality in \eqref{2.7} holds. Let $P_t(\mu,\nu,\d\gamma_1,\d\gamma_2)$ be the coupling of $Q_t(\mu,\d\gamma_1)$ and $Q_t(\nu,\d\gamma_2)$ introduced in the proof of Theorem~\ref{t2.2}. For any Borel function $F$ on $M(E)$ with $|F|\le 1$, we have
 \beqnn
\big|Q_tF(\nu)-Q_tF(\mu)\big|
 \ar=\ar
\bigg|\int_{M(E)^2} [F(\gamma_1) - F(\gamma_2)] P_t(\mu,\nu,\d\gamma_1,\d\gamma_2)\bigg| \cr
 \ar\le\ar
\int_{M(E)}Q_t(\mu\land\nu,\d\eta_0) \int_{M(E)}Q_t((\mu-\nu)_+,\d\eta_1) \cr
 \ar\ar\quad
\int_{M(E)}|F(\eta_0+\eta_1)-F(\eta_0+\eta_2)|Q_t((\mu-\nu)_-,\d\eta_2) \cr
 \ar\le\ar
2\int_{M(E)}Q_t(\mu\land\nu,\d\eta_0) \int_{M(E)}Q_t((\mu-\nu)_+,\d\eta_1) \cr
 \ar\ar\quad
\int_{M(E)}1_{\{\eta_1+\eta_2\neq 0\}}Q_t((\mu-\nu)_-,\d\eta_2) \cr
 \ar=\ar
2\int_{M(E)}Q_t(\mu\land\nu,\d\eta_0) \int_{M(E)}1_{\{\eta\neq 0\}}Q_t(|\mu-\nu|,\d\eta) \cr
 \ar=\ar
2\int_{M(E)}1_{\{\eta\neq 0\}}Q_t(|\mu-\nu|,\d\eta)
 =
2(1-\e^{-|\mu-\nu|(\bar{V}_t)}),
 \eeqnn
where the last equality follows by \eqref{2.6}. Then we have the second inequality in \eqref{2.7}. \qed

\bgcorollary\label{t2.9} Suppose that Condition~\ref{t2.6} holds. Then $L_{\var}(Q_tF)\le 2\|\bar{V}_t\|\|F\|$ for any $t>0$ and bounded Borel function $F$ on $M(E)$. \edcorollary

\bgcorollary\label{t2.10} Suppose that Condition~\ref{t2.6} holds. Then $\|Q_t(\mu,\cdot)-\delta_0\|_{\var} = 2(1-\e^{-\mu(\bar{V}_t)})\to 0$ as $t\to \infty$ for every $\mu\in M(E)$ if and only if $\lim_{t\to \infty} \bar{V}_t(x)= 0$ for every $x\in E$. \edcorollary

By Corollary~\ref{t2.9}, under Condition~\ref{t2.6} for any $t>0$ the operator $Q_t$ maps bounded Borel functions on $M(E)$ into functions continuous in the total variation distance. Then the semigroup $(Q_t)_{t\ge 0}$ has the so-called \textit{strong Feller property} in the total variation distance.

\section{Dawson--Watanabe superprocesses}

 \setcounter{equation}{0}

Let $\xi = (\Omega, \mcr{F}, \mcr{F}_t, \xi_t, \mbf{P}_x)$ be a Borel right process in $E$ with transition semigroup $(P_t)_{t\ge 0}$. Let $b\in B(E)$ and $c\in B(E)^+$. Let $\eta(x,\d y)$ be a bounded kernel on $E$ and $H(x,\d\nu)$ a $\sigma$-finite kernel from $E$ to $M(E)^\circ:= M(E)\setminus \{0\}$ satisfying
 \beqnn
\sup_{x\in E}\int_{M(E)^\circ}\big[\nu(1)\land\nu(1)^2 + \nu_x(1)\big]
H(x,\d\nu)< \infty,
 \eeqnn
where $\nu_x(\d y)$ denotes the restriction of $\nu(\d y)$ to $E\setminus\{x\}$. For $x\in E$ and $f\in B(E)^+$ write
 {\small\beqlb\label{3.1}
\phi(x,f)= b(x)f(x) + c(x)f(x)^2 - \eta(x,f) + \int_{M(E)^\circ} [\e^{-\nu(f)} - 1 + f(x)\nu(\{x\})] H(x,\d\nu).
 \eeqlb\!\!}
We can also rewrite \eqref{3.1} into
 \beqlb\label{3.2}
\phi(x,f)= b(x)f(x) + c(x)f(x)^2 - \gamma(x,f) + \int_{M(E)^\circ} [\e^{-\nu(f)} - 1 + \nu(f)] H(x,\d\nu),
 \eeqlb
where
 \beqnn
\gamma(x,\d y) = \eta(x,\d y) + \int_{M(E)^\circ}\nu_x(\d y) H(x,\d\nu).
 \eeqnn
By Proposition~2.20 in Li (2011), for every $f\in B(E)^+$ there is a unique locally bounded positive solution $(t,x)\mapsto V_tf(x)$ to the integral evolution equation
 \beqlb\label{3.3}
V_tf(x) = P_tf(x) - \int_0^t\d s\int_E \phi(y,V_sf) P_{t-s}(x,\d y), \quad
x\in E, t\ge 0.
 \eeqlb
By Theorem~5.12 in Li (2011), we can define a Borel right transition semigroup $(Q_t)_{t\ge 0}$ on $M(E)$ by \eqref{2.1}. If $X$ is a Markov process in $M(E)$ with transition semigroup $(Q_t)_{t\ge 0}$, we call it a \textit{Dawson--Watanabe superprocess} with \textit{spatial motion} $\xi$ and \textit{branching mechanism} $\phi$. For simplicity, we also called $X$ a \textit{$(\xi,\phi)$-superprocess}. By Theorem~2.27 of Li (2011), this process satisfies Condition~\ref{t2.1} with $(\pi_t)_{t\ge 0}$ defined by
 \beqlb\label{3.4}
\pi_tf(x) = P_tf(x) + \int_0^t P_{t-s}(\gamma-b) \pi_sf(x)\d s, \quad
t\ge 0,x\in E.
 \eeqlb

\bgexample\label{e2.1} In the special case where $E$ is a singleton, we can identify $M(E)$ with $[0,\infty)$. Let $\phi_*$ be a \textit{spatially independent} branching mechanism given by
 \beqlb\label{3.5}
\phi_*(z) = b_*z + c_*z^2 + \int_{(0,\infty)}\big(\e^{-zu}-1+zu\big)m_*(\d u), \quad
z\ge 0,
 \eeqlb
where $c_*\ge 0$ and $b_*$ are constants and $(u\land u^2)m_*(\d u)$ is a finite measure on $(0,\infty)$. We can define a transition semigroup $(Q_t)_{t\ge 0}$ by
 \beqlb\label{3.6}
\int_{[0,\infty)} \e^{-\lambda y} Q_t(x,\d y)
 =
\e^{-xv_t^*(\lambda)}, \qquad \lambda\ge 0,x\ge 0,
 \eeqlb
where $t\mapsto v_t^*(\lambda)$ is the unique positive solution of
 \beqlb\label{3.7}
\frac{\partial}{\partial t}v_t(\lambda)
 =
-\phi_*(v_t(\lambda)),
 \qquad
v_0(\lambda) = \lambda.
 \eeqlb
A Markov process $X$ in $[0,\infty)$ with transition semigroup $(Q_t)_{t\ge 0}$ is called a \textit{CB-process} with \textit{branching mechanism} $\phi$. The \textit{cumulant semigroup} of $X$ refers to the family of functions $(v_t)_{t\ge 0}$. From \eqref{3.6} and \eqref{3.7} it follows that
 \beqnn
\int_{[0,\infty)} y Q_t(x,\d y)= x\e^{-b_*t}, \qquad t\ge 0,x\ge 0,
 \eeqnn
which can be thought as a special form of \eqref{2.4}. We say the branching mechanism $\phi_*$ given by \eqref{3.5} satisfies \textit{Grey's condition} if $\phi_*(z)>0$ for sufficiently large $z>0$ and
 \beqlb\label{3.8}
\int^{\infty} \phi_*(z)^{-1}\d z< \infty.
 \eeqlb
For systematic studies of CB-processes, the reader may refer to Kyprianou (2014) and Li (2011, 2019+).
\edexample

\bgtheorem\label{t3.1} Let $(\pi_t)_{t\ge 0}$ be defined by \eqref{3.4} and let $\beta_*= \inf_{x\in E}[b(x) - \gamma(x,1)]$. Then $\|\pi_tf\|\le \e^{-\beta_*t}\|f\|$ for $t\ge 0$ and $f\in B(E)^+$.
\edtheorem

\proof Let $\beta(x)= b(x)-\gamma(x,1)$ for $x\in E$ and let $(P^\gamma_t)_{t\ge 0}$ be the locally bounded semigroup of kernels given by the Feynman--Kac formula
 \beqnn
P_t^\gamma f(x)=\mbf{P}_x\big[\e^{-\int_0^t\gamma(\xi_s,1)\d s} f(\xi_t)\big], \quad t\ge 0, x\in E, f\in B(E).
 \eeqnn
By Proposition~2.9 in Li (2011) we see \eqref{3.4} is equivalent to
 \beqlb\label{3.9}
\pi_tf(x) = P^\gamma_tf(x) + \int_0^t\d s\int_E [\gamma(y,\pi_sf) - \beta(y)\pi_sf(y)] P^\gamma_{t-s}(x,\d y).
 \eeqlb
By Theorem~A.43 in Li (2011) we can define a Borel right semigroup $(\tilde{P}_t)_{t\ge 0}$ on $E$ by
 \beqlb\label{3.10}
\tilde{P}_tf(x) = P^\gamma_tf(x) + \int_0^t\d s\int_E \gamma(y,\tilde{P}_sf) P^\gamma_{t-s}(x,\d y).
 \eeqlb
From \eqref{3.9} and \eqref{3.10} it follows that
 \beqnn
\pi_tf(x)= \tilde{P}_tf(x) - \int_0^t P^\gamma_{t-s}(\beta\pi_sf)(x)\d s + \int_0^t P^\gamma_{t-s}\gamma (\pi_sf-\tilde{P}_sf)(x)\d s.
 \eeqnn
Using the above relation successively we have
 \beqlb\label{3.11}
\pi_tf(x) \ar=\ar \tilde{P}_tf(x) - \int_0^t P^\gamma_{t-s_1} (\beta\pi_{s_1}f)(x)\d s_1 - \int_0^t\d s_1\int_0^{s_1} P^\gamma_{t-s_1}\gamma P^\gamma_{s_1-s_2}(\beta\pi_{s_2}f)(x)\d s_2 \cr
 \ar\ar
+ \int_0^t\d s_1\int_0^{s_1} P^\gamma_{t-s_1}\gamma P^\gamma_{s_1-s_2} \gamma(\pi_{s_2}f - \tilde{P}_{s_2}f)(x)\d s_2 \cr
 \ar=\ar
\tilde{P}_tf(x) - \int_0^t P^\gamma_{t-s_1} (\beta\pi_{s_1}f)(x)\d s_1 - \int_0^t\d s_1\int_0^{s_1} P^\gamma_{t-s_1}\gamma P^\gamma_{s_1-s_2} (\beta\pi_{s_2}f)(x)\d s_2 \cr
 \ar\ar
- \sum_{i=3}^n\int_0^t\d s_1\int_0^{s_1}\cdots \int_0^{s_{i-1}} P^\gamma_{t-s_1}\gamma P^\gamma_{s_1-s_2}\cdots \gamma P^\gamma_{s_{i-1}-s_i} (\beta\pi_{s_i}f)(x)\d s_i \ccr
 \ar\ar
+\, \varepsilon_n(t,x),
 \eeqlb
where
 \beqnn
\varepsilon_n(t,x)= \int_0^t\d s_1\int_0^{s_1}\cdots \int_0^{s_{n-1}} P^\gamma_{t-s_1}\gamma P^\gamma_{s_1-s_2}\cdots \gamma P^\gamma_{s_{n-1}-s_n} \gamma(\pi_{s_n}f - \tilde{P}_{s_n}f)(x)\d s_n.
 \eeqnn
By Proposition A.49 in Li (2011), there is a constant $a\ge 0$ so that $\|\pi_tf\|\le \|f\|\e^{at}$. Since $\|\tilde{P}_tf\|\le \|f\|$, we get
 \beqnn
\|\varepsilon_n(t,\cdot)\|
 \ar\le\ar
(1+\e^{at})\|f\|\|\gamma(\cdot,1)\|^n\int_0^t\d s_1\int_0^{s_1}\d s_2
\cdots\int_0^{s_{n-1}}\d s_n \cr
 \ar\le\ar
(1+\e^{at})\|f\|\|\gamma(\cdot,1)\|^n\frac{t^n}{n!}.
 \eeqnn
By Proposition~A.41 in Li (2011), the unique solution of \eqref{3.10} is given by
 \beqlb\label{3.12}
\tilde{P}_tf(x)
 =
P^\gamma_tf(x) + \sum_{i=1}^\infty\int_0^t\d s_1\int_0^{s_1}\d s_2\cdots \int_0^{s_{i-1}} P^\gamma_{t-s_1}\gamma P^\gamma_{s_1-s_2}\cdots \gamma P^\gamma_{s_i}f(x)\d s_i.
 \eeqlb
Then letting $n\to \infty$ in \eqref{3.11} and using \eqref{3.12} we obtain
 \beqlb\label{3.13}
\pi_tf(x) = \tilde{P}_tf(x) - \int_0^t\d s\int_E \beta(y)\pi_sf(y) \tilde{P}_{t-s}(x,\d y).
 \eeqlb
Let $\tilde{\xi} = (\tilde{\Omega}, \tilde{\mcr{F}}, \tilde{\mcr{F}}_t, \tilde{\xi}_t, \tilde{\mbf{P}}_x)$ be a right process realization of the Borel right semigroup $(\tilde{P}_t)_{t\ge 0}$. In view of \eqref{3.13}, we have
 \beqnn
\pi_tf(x)= \tilde{\mbf{P}}_x\big[\e^{-\int_0^t\beta(\tilde{\xi}_s)\d s} f(\tilde{\xi}_t)\big], \quad t\ge 0, x\in E.
 \eeqnn
Then $\|\pi_tf\|\le \e^{-\beta_*t}\tilde{\mbf{P}}_x[f(\tilde{\xi}_t)]\le\e^{-\beta_*t}\|f\|$ for $t\ge 0$. \qed

\bgcorollary\label{t3.2} Let $(Q_t)_{t\ge 0}$ be the transition semigroup of the $(\xi,\phi)$-superprocess defined by \eqref{2.1} and \eqref{3.3}. Then $W_1(Q_t(\mu,\cdot),\delta_0)= \mu(\pi_t1)\le \e^{-\beta_*t}\mu(1)$ for $t\ge 0$ and $\mu\in M(E)$.
\edcorollary

By Corollary~\ref{t3.2}, if $\beta_*> 0$, the transition law $Q_t(\mu,\cdot)$ converges to the stationary distribution $\delta_0$ exponentially fast in the Wasserstein distance as $t\to \infty$.

We next discuss the ergodicity in the total variation distance. The \textit{local projection} of the branching mechanism $\phi$ given by \eqref{3.1} or {\eqref{3.2} is the function $\phi_1$ on $E\times [0,\infty)$ defined by
 \beqlb\label{3.14}
\phi_1(x,z)= [b(x)-\gamma(x,1)]z + c(x)z^2 + \int_{M(E)^\circ} [\e^{-z\nu(\{x\})} - 1 + z\nu(\{x\})] H(x,\d\nu).
 \eeqlb
We say the branching mechanism $\phi$ is \textit{local} if $\gamma(\cdot,1)\equiv 0$. In this case, we also call $\phi_1$ the \textit{branching mechanism} of $X$.

\bgcondition\label{t3.3} The local projection $\phi_1$ of the branching mechanism is bounded below by a branching mechanism $\phi_*$ in the form \eqref{3.5}, that is, we have $\phi_1(x,z)\ge \phi_*(z)$ for all $x\in E$ and $z\ge 0$. \edcondition

\bgtheorem\label{t3.4} Suppose that Condition~\ref{t3.3} holds. Let $(V_t)_{t\ge 0}$ and $(v^*_t)_{t\ge 0}$ be defined by \eqref{3.3} and \eqref{3.7}, respectively. Then $\|V_tf\|\le v^*_t(\|f\|)$ for $t\ge 0$ and $f\in B(E)^+$. \edtheorem

\proof Up to an extension of the space $E$ as in the proof of Theorem 5.12 of Li (2011), we can assume $(P_t)_{t\ge 0}$ is a conservative transition semigroup. Let $\tilde{\phi}$ be the branching mechanism defined by
 \beqnn
\tilde{\phi}(x,f)= \phi_1(x,f(x)) + \gamma(x,1)f(x) - \gamma(x,f).
 \eeqnn
Let $(\tilde{V}_t)_{t\ge 0}$ denote the cumulant semigroup of the $(\xi,\tilde{\phi})$-superprocess. Then $(t,x)\mapsto \tilde{V}_tf(x)$ is the unique locally bounded positive solution to
 \beqlb\label{3.15}
\tilde{V}_tf(x) = P_tf(x) - \int_0^t\d s\int_E \tilde{\phi}(y,\tilde{V}_sf) P_{t-s}(x,\d y), \quad t\ge 0, x\in E.
 \eeqlb
It is easy to see that $\phi(x,f)\ge \tilde{\phi}(x,f)$ for $x\in E$ and $f\in B(E)^+$. By Corollary~5.18 in Li (2011) we have $V_tf(x)\le \tilde{V}_t f(x)$ for $x\in E$ and $f\in B(E)^+$. Let $(P^\gamma_t)_{t\ge 0}$ and $(\tilde{P}_t)_{t\ge 0}$ be the semigroups defined as in the proof of Theorem~\ref{t3.1}. By Proposition~2.9 in Li (2011), we can rewrite \eqref{3.15} into
 \beqlb\label{3.16}
\tilde{V}_tf(x) = P^\gamma_tf(x) - \int_0^t\d s\int_E \big[\phi_1 (y,\tilde{V}_sf(y)) - \gamma(y,\tilde{V}_sf)\big] P^\gamma_{t-s}(x,\d y).
 \eeqlb
From \eqref{3.10} and \eqref{3.16} it follows that
 \beqnn
\tilde{V}_tf(x) = \tilde{P}_tf(x) - \int_0^t P^\gamma_{t-s} \phi_1(\tilde{V}_sf)(x)\d s + \int_0^t P^\gamma_{t-s}\gamma (\tilde{V}_sf - \tilde{P}_sf)(x)\d s.
 \eeqnn
Using the above relation successively and arguing as in the proof of Theorem~\ref{t3.1} we see $(t,x)\mapsto \tilde{V}_tf(x)$ is also the unique locally bounded positive solution to
 \beqnn
\tilde{V}_tf(x) = \tilde{P}_tf(x) - \int_0^t\d s\int_E \phi_1(y,\tilde{V}_sf(y)) \tilde{P}_{t-s}(x,\d y).
 \eeqnn
Therefore we may think of $(\tilde{V}_t)_{t\ge 0}$ as the cumulant semigroup of a Dawson--Watanabe superprocess with local branching mechanism $\phi_1$ and underlying transition semigroup $(\tilde{P}_t)_{t\ge 0}$. Since $\phi_1(x,z)\ge \phi_*(z)$ for all $x\in E$ and $z\ge 0$, using Corollary 5.18 in Li (2011) again we see $\tilde{V}_tf(x)\le \tilde{V}_t\|f\|(x)\le v^*_t(\|f\|)$ for $t\ge 0$, $x\in E$ and $f\in B(E)^+$. \qed

\bgcorollary\label{t3.5} Suppose that Condition~\ref{t3.3} holds and $\phi_*^\prime(z)\to \infty$ as $z\to \infty$. Then we have \eqref{2.3} with $\lambda_t(x,1)= 0$ for $t>0$ and $x\in E$.
\edcorollary

\proof By Theorem~\ref{t3.4} we have $V_tf(x)\le \tilde{V}_tf(x)\le v^*_t(\|f\|)$ for $x\in E$ and $f\in B(E)^+$. Then the result follows as in the proof of Theorem~8.6 in Li (2011). \qed

\bgcorollary\label{t3.6} Suppose that Condition~\ref{t3.3} holds with $\phi_*$ satisfying Grey¡¯s condition \eqref{3.8}. Then Condition~\ref{t2.6} holds with $\|\bar{V}_t\|\le \bar{v}^*_t:= \lim_{\lambda\to \infty}v^*_t(\lambda)< \infty$ for $t>0$.
\edcorollary

\proof Since $\phi_*$ satisfies Grey's condition \eqref{3.11}, by Theorem~3.7 in Li (2011) we have $\bar{v}^*_t:= \lim_{\lambda\to \infty}v^*_t(\lambda)< \infty$ for $t>0$, so Theorem~\ref{t3.4} implies $\|\bar{V}_t\|\le \bar{v}^*_t< \infty$. \qed

\bgcorollary\label{t3.7} Suppose that Condition~\ref{t3.3} holds with $\phi_*$
satisfying Grey¡¯s condition \eqref{3.8}. Let $(Q_t)_{t\ge 0}$ be the transition semigroup of the $(\xi,\phi)$-superprocess defined by \eqref{2.1} and \eqref{3.3}. If $\beta_*:= \inf_{x\in E}[b(x) - \gamma(x,1)]>0$, then there is a constant $C\ge 0$ so that $\|Q_t(\mu,\cdot) - \delta_0\|_{\var}\le C(1+\mu(1))\e^{-\beta_*t}$ for $t\ge 0$ and $\mu\in M(E)$.
\edcorollary

\proof For $0\le t\le 1$ we have $\|Q_t(\mu,\cdot) - \delta_0\|_{\var}\le 2\le 2\e^{\beta_*}\e^{-\beta_*t}$. For any $t\ge 1$ we can use Theorem~\ref{t3.1} and Corollary~\ref{t3.6} to see $\|\bar{V}_t\|= \|V_{t-1}\bar{V}_1\|\le \|\pi_{t-1}\bar{V}_1\|\le \|\bar{V}_1\| \e^{-\beta_*(t-1)}$, and so $\|Q_t(\mu,\cdot) - \delta_0\|_{\var}\le 2\e^{\beta_*} \|\bar{V}_1\| \e^{-\beta_*t}\mu(1)$. Then we get the desired estimate with $C=2\e^{\beta_*}(1\vee \|\bar{V}_1\|)$. \qed

\bgcorollary\label{t3.8} Suppose that $\beta_*= \inf_{x\in E}[b(x) - \gamma(x,1)]> 0$ and $c_*= \inf_{x\in E}c(x)> 0$. Then we have
 \beqnn
\|V_tf\|\le \frac{\e^{-\beta_*t}\|f\|}{1 + c_*q(\beta_*,t)\|f\|},
 \qquad
t\ge 0, f\in B(E)^+,
 \eeqnn
where $q(\beta_*,t)= \beta_*^{-1}(1-\e^{-\beta_*t})$ with $q(0,t)= t$ by convention.
\edcorollary

\proof It is easy to see that Condition~\ref{t3.3} holds with $\phi_*(z) = \beta_*z + c_*z^2$. In this case, the solution of \eqref{3.7} is given by
 \beqnn
v_t^*(\lambda)
 =
\frac{\e^{-\beta_*t}\lambda}{1 + c_*q(\beta_*,t) \lambda}, \qquad t\ge 0, \lambda\ge 0.
 \eeqnn
Then the result follows by Theorem~\ref{t3.4}. \qed

Clearly, under the conditions of Corollary~\ref{t3.8}, we have $\|\bar{V}_t\|\le c_*^{-1} q(\beta_*,t)^{-1}\e^{-\beta_*t}$ for every $t>0$.

\section{MBI-processes and ergodicities}

 \setcounter{equation}{0}

Let $X$ be a MB-process with transition semigroup $(Q_t)_{t\ge 0}$. A generalization of the model can be formulated by introducing an immigration structure. A family of probability measures $(N_t)_{t\ge 0}$ on $M(E)$ is called a \textit{skew convolution semigroup} (SC-semigroup) associated with $X$ or $(Q_t)_{t\ge 0}$ provided
 \beqlb\label{4.1}
N_{r+t} = (N_rQ_t)*N_t, \qquad r,t\ge 0.
 \eeqlb
By Theorem~9.1 in Li (2011), the above relation is satisfied if and only if we can define another transition semigroup $(Q^N_t)_{t\ge 0}$ on $M(E)$ by
 \beqlb\label{4.2}
Q^N_t(\mu,\cdot) = Q_t(\mu,\cdot)*N_t, \qquad t\ge 0,\ \mu\in M(E).
 \eeqlb
A Markov process in $M(E)$ with transition semigroup $(Q^N_t)_{t\ge 0}$ is naturally called an \textit{MBI-process} associated with $X$ or $(Q_t)_{t\ge 0}$.

By Theorem~9.4 of Li (2011), there is a one-to-one correspondence between SC-semigroups $(N_t)_{t\ge 0}$ and infinitely divisible probability entrance laws $(K_t)_{t>0}$ for the semigroup $(Q_t)_{t\ge 0}$ satisfying
 \beqnn
-\int_0^t \log L_{K_s}(1)\d s<\infty, \qquad t\ge 0.
 \eeqnn
The one-to-one correspondence is determined by
 \beqlb\label{4.3}
L_{N_t}(f)= \exp\bigg\{-\int_0^t I_s(K,f)\d s\bigg\}, \qquad t\ge 0,f\in B(E)^+,
 \eeqlb
where $I_s(K,f) = - \log L_{K_s}(f)$. For the SC-semigroup $(N_t)_{t\ge 0}$ represented by \eqref{4.3}, the corresponding transition semigroup $(Q_t^N)_{t\ge 0}$ defined in \eqref{4.2} is given by
 \beqlb\label{4.4}
\int_{M(E)}\e^{-\nu(f)}Q_t^N(\mu,\d\nu)
 =
\exp\bigg\{-\mu(V_tf) -\int_0^tI_s(K,f)\d s\bigg\}.
 \eeqlb

We are particularly interested in SC-semigroups with finite first moments. By replacing $f\in B(E)^+$ with $\lambda f$ for $\lambda\ge 0$ in \eqref{4.3} and taking the right derivatives at $\lambda=0$ we get
 \beqlb\label{4.5}
\int_{M(E)}\nu(f)N_t(\d\nu) = \int_0^t\d s\int_{M(E)}\nu(f)K_s(\d\nu).
 \eeqlb
Then $N_t$ has finite first moment if and only if
 \beqlb\label{4.6}
\int_0^t\d s\int_{M(E)}\nu(1)K_s(\d\nu)< \infty.
 \eeqlb
In this case, we can extend \eqref{4.5} to all $f\in B(E)$. If Condition~\ref{t2.1} also holds, we have
 \beqlb\label{4.7}
\int_{M(E)}\nu(f)Q^N_t(\mu,\d\nu) = \mu(\pi_tf) + \int_0^t\d s\int_{M(E)}\nu(f)K_s(\d\nu).
 \eeqlb

\bgtheorem\label{t4.1} Suppose that \eqref{4.6} and Condition~\ref{t2.1} hold. Let $(Q^N_t)_{t\ge 0}$ be the transition semigroup defined by \eqref{4.4}. Then, for $t\ge 0$ and $\mu,\nu\in M(E)$,
\beqlb\label{4.8}
|\mu(\pi_t1)-\nu(\pi_t1)|
 \le
W_1(Q_t^N(\mu,\cdot),Q_t^N(\nu,\cdot))
 \le
|\mu-\nu|(\pi_t1).
 \eeqlb
\edtheorem

\proof The first inequality in \eqref{4.8} follows by a first moment calculation based on \eqref{4.7} as in the proof of Theorem~\ref{t2.2}. Let $P_t(\mu,\nu,\d\gamma_1,\d\gamma_2)$ be the coupling of $Q_t(\mu,\d\gamma_1)$ and $Q_t(\nu,\d\gamma_2)$ defined in that proof. Let $Q_t(\mu,\nu, \d\eta_1,\d\eta_2)$ be the image of $N_t(\d\gamma_0) P_t(\mu,\nu, \d\gamma_1,\d\gamma_2)$ under the mapping $(\gamma_0,\gamma_1,\gamma_2)\mapsto (\eta_1,\eta_2):= (\gamma_0+\gamma_1, \gamma_0+\gamma_2)$. From \eqref{4.2} we see that $Q_t(\mu,\nu, \d\eta_1,\d\eta_2)$ is a coupling of $Q_t^N(\mu,\d\eta_1)$ and $Q_t^N(\nu,\d\eta_2)$. It follows that
 \beqnn
W_1(Q_t^N(\mu,\cdot),Q_t^N(\nu,\cdot))
 \ar\le\ar
\int_{M(E)^2} \|\eta_1-\eta_2\|_{\var} Q_t(\mu,\nu,\d\eta_1,\d\eta_2) \cr
 \ar=\ar
\int_{M(E)} N_t(\d\gamma_0)\int_{M(E)^2} \|\gamma_1-\gamma_2\|_{\var} P_t(\mu,\nu,\d\gamma_1,\d\gamma_2) \cr
 \ar=\ar
\int_{M(E)^2} \|\gamma_1-\gamma_2\|_{\var} P_t(\mu,\nu,\d\gamma_1,\d\gamma_2).
 \eeqnn
Then the second inequality in \eqref{4.8} follows by the calculations in the proof of Theorem~\ref{t2.2}. \qed

\bgcorollary\label{t4.2} Suppose that \eqref{4.6} and Condition~\ref{t2.1} hold. Let $F$ be a Borel function on $M(E)$. Then $L_{\var}(Q_t^NF)\le \|\pi_t1\|L_{\var}(F)$ for $t\ge 0$. \edcorollary

\bgtheorem\label{t4.3} Suppose that Condition~\ref{t2.6} holds. Let $(Q^N_t)_{t\ge 0}$ be the transition semigroup defined by \eqref{4.2}. Then, for $t>0$ and $\mu,\nu\in M(E)$,
 \beqnn
\|Q_t^N(\mu,\cdot)-Q_t^N(\nu,\cdot)\|_{\var}
 \le
2(1-\e^{-|\mu-\nu|(\bar{V}_t)})
 \le
2|\mu-\nu|(\bar{V}_t).
 \eeqnn
\edtheorem

\proof Let $F$ be a Borel function on $M(E)$ satisfying $|F|\le 1$. In view of \eqref{4.2}, we have
 \beqnn
\big|Q_t^NF(\mu)-Q_t^NF(\nu)\big|
 \ar=\ar
\bigg|\int_{M(E)} F(\eta)Q_t^N(\mu,\d\eta) - \int_{M(E)} F(\eta)Q_t^N(\nu,\d\eta)\bigg| \cr
 \ar=\ar
\bigg|\int_{M(E)} N_t(\d\gamma)\int_{M(E)} F(\eta+\gamma)Q_t(\mu,\d\eta) \cr
 \ar\ar\qquad\qquad
- \int_{M(E)} N_t(\d\gamma)\int_{M(E)} F(\eta+\gamma)Q_t(\nu,\d\eta)\bigg| \cr
 \ar\le\ar
\int_{M(E)}\bigg|\int_{M(E)} F(\eta+\gamma)Q_t(\mu,\d\eta) \cr
 \ar\ar\qquad\qquad
- \int_{M(E)} F(\eta+\gamma)Q_t(\nu,\d\eta)\bigg|N_t(\d\gamma) \ccr
 \ar\le\ar
\big\|Q_t(\nu,\cdot) - Q_t(\nu,\cdot)\big\|_{\var}.
 \eeqnn
Then desired estimates follow by Theorem~\ref{t2.8}. \qed

\bgcorollary\label{t4.4} Suppose that Condition~\ref{t2.6} holds. Let $F$ be a Borel function on $M(E)$. Then $L_{\var}(Q_t^NF)\le 2\|\bar{V}_t\|\|F\|$ for $t>0$. \edcorollary

By Corollary~\ref{t4.4}, under Condition~\ref{t2.6} the semigroup $(Q_t^N)_{t\ge 0}$ is strong Feller in the total variation distance.

\bgtheorem\label{t4.5} Let $(N_t)_{t\ge 0}$ be the SC-semigroup given by \eqref{4.3}. Then $N_t$ converges weakly as $t\to \infty$ to a probability measure $N_\infty$ on $M(E)$ with finite first moment if and only if
 \beqlb\label{4.9}
\int_0^\infty\d s\int_{M(E)}\nu(1)K_s(\d\nu)< \infty.
 \eeqlb
In this case, we have
 \beqlb\label{4.10}
L_{N_\infty}(f)= \exp\bigg\{-\int_0^\infty I_s(K,f)\d s\bigg\}, \qquad f\in B(E)^+
 \eeqlb
and
 \beqlb\label{4.11}
\int_{M(E)^\circ} \nu(f)N_\infty(\d\nu)
 =
\int_0^\infty \d s\int_{M(E)}\nu(f)K_s(\d\nu), \qquad f\in B(E).
 \eeqlb
\edtheorem

\proof Suppose that \eqref{4.9} holds. By Jensen's inequality, for $f\in B(E)^+$ we have
 \beqnn
-\lim_{t\to \infty}\int_0^t \log L_{K_s}(f)\d s
 \ar=\ar
-\int_0^\infty \log L_{K_s}(f)\d s \cr
 \ar\le\ar
\int_0^\infty\d s\int_{M(E)}\nu(f)K_s(\d\nu)< \infty.
 \eeqnn
Clearly, the convergence above is uniform on $\{f\in B(E)^+: \|f\|\le a\}$ for each $a\ge 1$. By Corollary~1.21 in Li (2011), we can define a probability measure $N_\infty$ by \eqref{4.10} and $\lim_{t\to \infty}N_t= N_\infty$ by weak convergence. Then we get \eqref{4.11} from \eqref{4.10}. Conversely, suppose that $N_t$ converges weakly as $t\to \infty$ to a probability measure $N_\infty$ on $M(E)$ with finite first moment. By \eqref{4.3} we see \eqref{4.10} holds for continuous functions $f\in B(E)^+$, so it holds all $f\in B(E)^+$. From \eqref{4.10} we get \eqref{4.11}. Then \eqref{4.9} is satisfied. \qed

\bgcorollary\label{t4.6} Suppose that \eqref{4.9} holds. Then the probability $N_\infty$ defined by \eqref{4.10} is a stationary distribution for the semigroup $(Q^N_t)_{t\ge 0}$.
\edcorollary

\bgcorollary\label{t4.7} Suppose that \eqref{4.9} holds. Then $\lim_{t\to \infty}Q^N_t(\mu,\cdot)= N_\infty$ by weak convergence for every $\mu\in M(E)$ if and only if $\lim_{t\to \infty}V_t1(x)= 0$ for every $x\in E$. In this case, $N_\infty$ is the unique stationary distribution for $(Q^N_t)_{t\ge 0}$. \edcorollary

The next theorem gives an accurate evaluations of the distances between the SC-semigroup $(N_t)_{t\ge 0}$ and its limit distribution:

\bgtheorem\label{t4.8} Let $(N_t)_{t\ge 0}$ be the SC-semigroup given by \eqref{4.3}. Suppose that \eqref{4.9} and Condition~\ref{t2.1} hold. Then for $t\ge 0$ we have
 \beqnn
W_1(N_t,N_\infty)
 =
\int_t^\infty\d s \int_{M(E)} \nu(1) K_s(\d\nu)
 =
\int_{M(E)}\nu(\pi_t1) N_\infty(\d\nu).
 \eeqnn
\edtheorem

\proof Since $(K_t)_{t>0}$ is an entrance law for $(Q_t)_{t\ge 0}$, the second desired equality follows from \eqref{2.4}. By \eqref{4.3} and \eqref{4.10} it is easy to show that $N_t*(N_\infty Q_t)= N_\infty$ for $t\ge 0$. Let $M_t(\d\eta_1,\d\eta_2)$ be the image of the product measure $N_t(\d\nu_1)(N_\infty Q_t)(\d\nu_2)$ under the mapping $(\nu_1,\nu_2)\mapsto (\eta_1,\eta_2):= (\nu_1, \nu_1+\nu_2)$. Then $M_t(\d\eta_1,\d\eta_2)$ is a coupling of $N_t(\d\eta_1)$ and $N_\infty(\d\nu_2)$. It follows that
 \beqnn
W_1(N_t,N_\infty)
 \ar\le\ar
\int_{M(E)^2} \|\eta_1-\eta_2\|_{\var} M_t(\d\eta_1,\d\eta_2) \cr
 \ar=\ar
\int_{M(E)} N_t(\d\nu_1)\int_{M(E)} \|\nu_2\|_{\var} N_\infty Q_t(\d\nu_2) \cr
 \ar=\ar
\int_{M(E)} N_\infty(\d\nu)\int_{M(E)} \nu_2(1) Q_t(\nu,\d\nu_2) \cr
 \ar=\ar
\int_0^\infty\d s \int_{M(E)} K_s(\d\nu) \int_{M(E)} \nu_2(1)Q_t(\nu,\d\nu_2) \cr
 \ar=\ar
\int_0^\infty\d s \int_{M(E)} \nu(1)K_{s+t}(\d\nu) \cr
 \ar=\ar
\int_t^\infty\d s \int_{M(E)} \nu(1)K_s(\d\nu),
 \eeqnn
where we have used \eqref{4.11} for the third equality. On the other hand, by Theorem~5.10 in Chen (2004a, p.181) we have
 \beqnn
W_1(N_t,N_\infty)
 \ar\ge\ar
\int_{M(E)} \eta(1) (N_\infty-N_t)(\d\eta) \cr
 \ar=\ar
\int_0^\infty\d s \int_{M(E)} \eta(1)K_s(\d\eta) - \int_0^t\d s \int_{M(E)} \eta(1)K_s(\d\eta) \cr
 \ar=\ar
\int_t^\infty\d s \int_{M(E)} \nu(1)K_s(\d\nu),
 \eeqnn
where we have also used \eqref{4.5} and \eqref{4.11}. Then the desired relation holds. \qed

\bgcorollary\label{t4.9} Suppose that \eqref{4.9} and Condition~\ref{t2.1} hold. Then for $t\ge 0$ and $\mu\in M(E)$ we have
 \beqnn
W_1(Q_t^N(\mu,\cdot),N_\infty)
 \le
\mu(\pi_t1) + \int_{M(E)}\nu(\pi_t1) N_\infty(\d\nu).
 \eeqnn
\edcorollary

\proof By \eqref{4.2} we have $N_t=Q_t^N(0,\cdot)$. Then the result follows by the estimates given in  Theorems~\ref{t4.1} and~\ref{t4.8} and the triangle inequality. \qed

\bgtheorem\label{t4.10} Suppose that \eqref{4.9} and Condition~\ref{t2.6} hold. Then for $t>0$ we have
 \beqnn
\|N_t-N_\infty\|_{\var}
 \le
2\int_{M(E)}(1-\e^{-\nu(\bar{V}_t)})N_\infty(\d\nu)
 \le
2\int_{M(E)}\nu(\bar{V}_t) N_\infty(\d\nu).
 \eeqnn
\edtheorem

\proof Let $M_t(\d\eta_1,\d\eta_2)$ be the coupling of $N_t(\d\eta_1)$ and $N_\infty(\d\nu_2)$ defined in the proof of Theorem~\ref{t4.8}. By Theorem~5.7 in Chen (2004a, p.179) we have
 \beqnn
\|N_t-N_\infty\|_{\var}
 \ar\le\ar
2\int_{M(E)^2} 1_{\{\eta_1\neq \eta_2\}} M_t(\d\eta_1,\d\eta_2) \cr
 \ar=\ar
2\int_{M(E)} N_t(\d\nu_1)\int_{M(E)} 1_{\{\nu_2\neq 0\}} N_\infty Q_t(\d\nu_2) \cr
 \ar=\ar
2\int_{M(E)} N_\infty(\d\nu)\int_{M(E)} 1_{\{\nu_2\neq 0\}} Q_t(\nu,\d\nu_2) \cr
 \ar=\ar
2\int_{M(E)}(1-\e^{-\nu(\bar{V}_t)})N_\infty(\d\nu) \cr
 \ar\le\ar
2\int_{M(E)} \nu(\bar{V}_t) N_\infty(\d\nu),
 \eeqnn
where we used \eqref{2.6} for the last equality. \qed

\bgcorollary\label{t4.11} Suppose that \eqref{4.9} and Condition~\ref{t2.6} hold. Then, for $t>0$ and $\mu\in M(E)$,
 \beqnn
\|Q_t^N(\mu,\cdot)-N_\infty\|_{\var}
 \le
2\mu(\bar{V}_t) + 2\int_{M(E)}\nu(\bar{V}_t) N_\infty(\d\nu).
 \eeqnn
\edcorollary

\proof This follows by Theorems~\ref{t4.3} and~\ref{t4.10} and the triangle inequality. \qed

In the sequel of this section, we consider the $(\xi,\phi)$-superprocess with the transition semigroup $(Q_t)_{t\ge 0}$ defined by \eqref{2.1} and \eqref{3.3}. Let $\mcr{K}(P)$ and $\mcr{K}(\pi)$ denote the set of entrance laws $\kappa= (\kappa_t)_{t>0}$ for the semigroups $(P_t)_{t\ge 0}$ and $(\pi_t)_{t\ge 0}$, respectively, satisfying the integrability condition
 \beqnn
\int_0^1\kappa_s(1)\d s< \infty.
 \eeqnn
Given $\kappa\in \mcr{K}(P)$ we set, for $t>0$ and $f\in B(E)^+$,
 \beqlb\label{4.12}
\pi_t(\kappa,f) = \kappa_t(f) + \int_0^t \kappa_{t-s}((\gamma-b)\pi_sf)\d s
 \eeqlb
and
 \beqlb\label{4.13}
V_t(\kappa,f) = \kappa_t(f) - \int_0^t\d s\int_E \phi(y,V_sf)\kappa_{t-s}(\d y).
 \eeqlb
Clearly, if $\kappa= (\kappa_t)_{t>0} $ is \textit{closed} by a measure $\mu$ on $E$ in the sense $\kappa_t= \mu P_t$, then $\pi_t(\kappa,f)= \mu(\pi_tf)$ and $V_t(\kappa,f)= \mu(V_tf)$. The reader may refer to Dynkin (1989) and Li (1996b, 2011) for the discussions on the connections between entrance laws for the $(\xi,\phi)$-superprocess and those for the underlying process.

\bgremark\label{t4.12} If $\kappa= (\kappa_t)_{t>0}\in \mcr{K}(P)$ or $\mcr{K}(\pi)$, then each $\kappa_t$ is a finite measure on $E$. Indeed, by the above integrability condition, for any $t>0$ we can find a sequence $r\in (0,t]$ so that $\kappa_r(1)< \infty$. In the case of $\kappa\in \mcr{K}(P)$, we have $\kappa_t(1) = \kappa_r(P_{t-r}1)\le \kappa_r(1)< \infty$. In the case of $\kappa\in \mcr{K}(\pi)$, by Theorem~\ref{t3.1} we have $\kappa_t(1) = \kappa_r(\pi_{t-r}1)\le \e^{-\beta_*(t-r)}\kappa_r(1)< \infty$. \edremark

We endow $\mcr{K}(P)$ with the $\sigma$-algebra generated by the collection of maps $\{\kappa\mapsto \kappa_t(f): t>0, f\in B(E)\}$. Given an entrance law $\kappa\in\mcr{K}(P)$ and a $\sigma$-finite measure $F(\d\nu)$ on $\mcr{K}(P)^\circ:= \mcr{K}(P)\setminus \{0\}$ satisfying
 \beqlb\label{4.14}
\int_0^1\d s\int_{\mcr{K}(P)^\circ} \nu_s(1) F(\d\nu)< \infty
 \eeqlb
write, for $t>0$ and $f\in B(E)^+$,
 \beqlb\label{4.15}
I_t(\kappa,F,f)= V_t(\kappa,f) + \int_{\mcr{K}(P)^\circ} (1-\e^{-V_t(\nu,f)}) F(\d\nu).
 \eeqlb

\bgtheorem\label{t4.13} There is a one-to-one correspondence between SC-semigroups $(N_t)_{t\ge 0}$ with finite first moments and the pairs $(\kappa,F)$, for $\kappa\in\mcr{K}(P)$ and for $F(\d\nu)$ satisfying \eqref{4.14}, given by
 \beqlb\label{4.16}
L_{N_t}(f)= \exp\bigg\{-\int_0^t I_s(\kappa,F,f)\d s\bigg\}, \qquad t\ge 0,f\in B(E)^+.
 \eeqlb
Moreover, if $(N_t)_{t\ge 0}$ and $(\kappa,F)$ are related by \eqref{4.16}, then, for $t\ge 0$ and $f\in B(E)^+$,
 \beqlb\label{4.17}
\int_{M(E)^\circ} \nu(f)N_t(\d\nu)
 =
\int_0^t \bigg[\pi_s(\kappa,f) + \int_{\mcr{K}(P)^\circ} \pi_s(\nu,f)F(\d\nu)\bigg] \d s.
 \eeqlb
\edtheorem

\proof Let $(N_t)_{t\ge 0}$ be given by \eqref{4.3}. By \eqref{2.4} and \eqref{4.5} one can see that $(N_t)_{t\ge 0}$ has finite first moments if and only if $(K_t)_{t>0}$ satisfies
 \beqnn
\int_0^1\d s\int_{M(E)}\nu(1)K_s(\d\nu)< \infty.
 \eeqnn
Let $\mcr{K}^1(Q)$ denote the set of probability entrance laws $(K_t)_{t>0}$ for $(Q_t)_{t\ge 0}$ satisfying the above integrability condition. By Theorem~8.20 of Li (2011), an infinitely divisible probability entrance laws $(K_t)_{t>0}\in \mcr{K}^1(Q)$ corresponds uniquely to a pair $(\kappa,F)$, where $\kappa\in\mcr{K}(P)$ and $F(\d\nu)$ satisfies \eqref{4.14}. The correspondence is given by
 \beqnn
L_{K_t}(f)= \exp\big\{-I_t(\kappa,F,f)\big\}, \qquad t>0,f\in B(E)^+.
 \eeqnn
Thus \eqref{4.16} establishes a one-to-one correspondence SC-semigroups $(N_t)_{t\ge 0}$ with finite first moments and the pairs $(\kappa,F)$. By (8.42) in Li (2011) we have
 \beqlb\label{4.18}
\int_{M(E)^\circ} \nu(f)K_t(\d\nu)
 =
\pi_t(\kappa,f) + \int_{\mcr{K}(P)^\circ} \pi_t(\nu,f)F(\d\nu).
 \eeqlb
Then \eqref{4.17} follows from \eqref{4.5}. \qed

For the SC-semigroup $(N_t)_{t\ge 0}$ represented by \eqref{4.16}, the corresponding transition semigroup $(Q_t^N)_{t\ge 0}$ defined in \eqref{4.2} is given by
 \beqlb\label{4.19}
\int_{M(E)}\e^{-\nu(f)}Q_t^N(\mu,\d\nu)
 =
\exp\bigg\{-\mu(V_tf) -\int_0^tI_s(\kappa,F,f)\d s\bigg\}.
 \eeqlb

We omit the proofs of some of the following results for the $(\xi,\phi)$-superprocess as they are easy consequences of the general results established in the first part of this section.

\bgtheorem\label{t4.14} Let $(N_t)_{t\ge 0}$ be the SC-semigroup defined by \eqref{4.16}. Then $N_t$ converges weakly as $t\to \infty$ to a probability measure $N_\infty$ on $M(E)$ with finite first moment if and only if
 \beqlb\label{4.20}
\int_0^\infty \bigg[\pi_s(\kappa,1) + \int_{\mcr{K}(P)^\circ} \pi_s(\nu,1)F(\d\nu)\bigg] \d s< \infty.
 \eeqlb
In this case, we have, for $f\in B(E)^+$,
 \beqlb\label{4.21}
L_{N_\infty}(f)
 =
\exp\bigg\{-\int_0^\infty I_s(\kappa,F,f)\d s\bigg\}
 \eeqlb
and
 \beqlb\label{4.22}
\int_{M(E)^\circ} \nu(f)N_\infty(\d\nu)
 =
\int_0^\infty \bigg[\pi_s(\kappa,f) + \int_{\mcr{K}(P)^\circ} \pi_s(\nu,f)F(\d\nu)\bigg] \d s.
 \eeqlb
\edtheorem

\bgtheorem\label{t4.15} Let $(N_t)_{t\ge 0}$ be the SC-semigroup defined by \eqref{4.16}. Suppose that \eqref{4.20} holds. Then for $t\ge 0$ we have
 \beqnn
W_1(N_t,N_\infty)
 =
\int_{M(E)}\nu(\pi_t1) N_\infty(\d\nu)
 =
\int_t^\infty \bigg[\pi_s(\kappa,1) + \int_{\mcr{K}(P)^\circ} \pi_s(\nu,1)F(\d\nu)\bigg] \d s.
 \eeqnn
\edtheorem

\bgcorollary\label{t4.16} Suppose that \eqref{4.20} holds. {\rm(i)}~We have $\lim_{t\to \infty} W_1(N_t,N_\infty)= 0$. {\rm(ii)}~We have $\lim_{t\to \infty} W_1(Q_t^N(\mu,\cdot),N_\infty)= 0$ for every $\mu\in M(E)$ if $\lim_{t\to \infty} \pi_t1(x)= 0$ for every $x\in E$. {\rm(iii)}~If $\beta_*:= \inf_{x\in E}[b(x) - \gamma(x,1)]>0$, there is a constant $C\ge 0$ so that $W_1(Q_t^N(\mu,\cdot),N_\infty)\le C(1+\mu(1)) \e^{-\beta_*t}$ for $t\ge 0$ and $\mu\in M(E)$.
\edcorollary

\proof By Theorem~\ref{t4.15} and Corollary~\ref{t4.9} we have (i) and (ii). The assertion (iii) follows from the estimate $\pi_t1(x)\le \e^{-\beta_*t}$ for $t\ge 0$ and $x\in E$. \qed

\bgtheorem\label{t4.17} Let $(N_t)_{t\ge 0}$ be the SC-semigroup defined by \eqref{4.16}. Suppose that \eqref{4.20} and Condition~\ref{t3.3} hold with $\phi_*$ satisfying Grey's condition \eqref{3.8}. Then for $t>0$ we have
 \beqnn
\|N_t-N_\infty\|_{\var}
 \le
2\int_{M(E)}\nu(\bar{V}_t) N_\infty(\d\nu)
 =
2\int_0^\infty \bigg[\pi_s(\kappa,\bar{V}_t) + \int_{\mcr{K}(P)^\circ} \pi_s(\nu,\bar{V}_t) F(\d\nu)\bigg] \d s.
 \eeqnn
\edtheorem

\bgcorollary\label{t4.18} Suppose that \eqref{4.20} and Condition~\ref{t3.3} hold with $\phi_*$ satisfying Grey's condition \eqref{3.8}. {\rm(i)}~We have $\lim_{t\to \infty}\|N_t-N_\infty\|_{\var}= 0$. {\rm(ii)}~If $\lim_{t\to \infty} \bar{V}_t(x)= 0$ for every $x\in E$, then $\lim_{t\to \infty}\|Q_t^N(\mu,\cdot) - N_\infty\|_{\var}= 0$ for every $\mu\in M(E)$. {\rm(iii)}~If $\beta_*:= \inf_{x\in E}[b(x) - \gamma(x,1)]>0$, then there is a constant $C\ge 0$ so that $\|Q_t^N(\mu,\cdot)-N_\infty\|_{\var}\le C(1+\mu(1))\e^{-\beta_*t}$ for $t\ge 0$ and $\mu\in M(E)$.
\edcorollary

\proof It is easy to see that $\bar{V}_t(x)= V_{t-r}\bar{V}_r(x)\le \pi_{t-r}\bar{V}_r(x)\le \|\bar{V}_r\| \pi_{t-r}1(x)$ for $t\ge r> 0$ and $x\in E$. From the estimate in Theorem~\ref{t4.17} it follows that
 \beqnn
\|N_t-N_\infty\|_{\var}
 \ar\le\ar
\|\bar{V}_r\|\int_0^\infty \bigg[\pi_s(\kappa,\pi_{t-r}1) + \int_{\mcr{K}(P)^\circ} \pi_s(\nu,\pi_{t-r}1)F(\d\nu)\bigg] \d s \cr
 \ar=\ar
\|\bar{V}_r\|\int_{t-r}^\infty \bigg[\pi_s(\kappa,1) + \int_{\mcr{K}(P)^\circ} \pi_s(\nu,1) F(\d\nu)\bigg] \d s.
 \eeqnn
The right-hand side goes to zero as $t\to \infty$. That gives (i). By Corollary~\ref{t4.11} we have (ii). The assertion (iii) follows as in the proof of Corollary~\ref{t3.7}. \qed

\bgremark\label{t4.19} Suppose that $U(x):= \int_0^\infty \pi_s1(x)\d s$ is bounded on $E$. Then \eqref{4.20} follows from our assumptions on the pair $(\kappa,F)$. Indeed, the quantity in \eqref{4.20} is equal to
 \beqnn
\int_0^1 \bigg[\pi_s(\kappa,1) + \int_{\mcr{K}(P)^\circ} \pi_s(\nu,1) F(\d\nu)\bigg] \d s + \pi_1(\kappa,U) + \int_{\mcr{K}(P)^\circ} \pi_s(\nu,U) F(\d\nu).
 \eeqnn
By Theorem~\ref{t3.1}, the function $U$ is bounded on $E$ if $\beta_*:= \inf_{x\in E}[b(x) - \gamma(x,1)]>0$.
\edremark

\bgremark\label{t4.20} Let $h\in B(E)^+$ be a strictly positive $\alpha$-excessive function for $(P_t)_{t\ge 0}$ for some $\alpha\ge 0$. We can define the transition semigroup $(\bar{P}_t)_{t\ge 0}$ of a Borel right process $\bar{\xi}$ on $E$ by
 \beqnn
\bar{P}_tf(x) = h(x)^{-1}\e^{-\alpha t}P_t(x,hf), \quad x\in E, f\in
B(E).
 \eeqnn
Starting from $\bar{\xi}$ as the underlying process, we can construct a Dawson--Watanabe type process $\bar{X}$ in the state space $M(E)$. Let $M_h(E)$ denote the space of tempered measures $\mu$ on $E$ satisfying $\mu(h)< \infty$. From $\bar{X}$ we can use the homeomorphic transformation $\mu(\d x) \mapsto h(x)^{-1} \mu(\d x)$ from $M(E)$ to $M_h(E)$ to obtain a Dawson--Watanabe type process $X$ in $M_h(E)$. The reader may refer to Section~6.1 of Li (2011) for the detailed arguments. By this transformation, the results obtained in this and the last two sections can be reformulated for the state space $M_h(E)$.
\edremark

\section{Self-decomposable distributions}

 \setcounter{equation}{0}

For probability measures $F$ and $H$ on $M(E)$, we write $F\preceq H$ if there is another probability measure $G$ on $M(E)$ so that $F*G = H$. Clearly, the probability $G$ is unique if it exists. Let $(Q_t)_{t\ge 0}$ be the transition semigroup of the $(\xi,\phi)$-superprocess. We say a probability $N$ on $M(E)$ is \textit{self-decomposable} or \textit{C-excessive} for $(Q_t)_{t\ge 0}$ if $NQ_t\preceq N$ for all $t\ge 0$. Let $\mcr{E}^*(Q)$ denote the set of C-excessive probabilities on $M(E)$. By Theorem~9.8 of Li (2011), for each $N\in \mcr{E}^*(Q)$ there is a unique SC-semigroup $(N_t)_{t\ge 0}$ associated with $(Q_t)_{t\ge 0}$ such that
 \beqlb\label{5.1}
N = (NQ_t)*N_t, \qquad t\ge 0.
 \eeqlb
It is easy to see that $\mcr{E}^*(Q)$ contains the set of stationary (or invariant) probabilities $\mcr{E}_i^*(Q)$ for $(Q_t)_{t\ge 0}$. We say $N\in \mcr{E}^*(Q)$ is \textit{purely self-decomposable} or \textit{purely C-excessive} if $\lim_{t\to \infty} NQ_t = \delta_0$ by weak convergence. Let $\mcr{E}_p^*(Q)$ denote the set of purely C-excessive probabilities for $(Q_t)_{t\ge 0}$. By Theorem~9.10 of Li (2011), a C-excessive probability $N\in \mcr{E}^*(Q)$ has the unique decomposition $N = N^i*N^p$, where $N^i= \lim_{t\to \infty} NQ_t\in \mcr{E}_i^*(Q)$ and $N^p = \lim_{t\to \infty}N_t\in \mcr{E}_p^*(Q)$. In particular, for $N\in \mcr{E}_p^*(Q)$ we have $N = \lim_{t\to \infty}N_t$. By Theorem~\ref{t4.14} there is a one-to-one correspondence given by \eqref{4.21} between the distributions $N=N_\infty\in \mcr{E}_p^*(Q)$ with finite first moment and the pairs $(\kappa,F)$ satisfying \eqref{4.20}.

A $\sigma$-finite measure $\gamma$ on $E$ is said to be \textit{excessive} for the semigroup $(\pi_t)_{t\ge 0}$ on $E$ if $\gamma\pi_t\le \gamma$ for all $t\ge 0$. Let $\mcr{E}(\pi)$ denote the set of all excessive finite measures for $(\pi_t)_{t\ge 0}$. We say $\gamma\in \mcr{E}(\pi)$ is \textit{purely excessive} if $\lim_{t\to \infty} \gamma Q_t = 0$. Let $\mcr{E}_p(\pi)\subset \mcr{E}(\pi)$ denote the set of purely excessive finite measures for $(\pi_t)_{t\ge 0}$. Let $\mcr{E}_i(\pi)\subset \mcr{E}(\pi)$ be the set of invariant finite measures for $(\pi_t)_{t\ge 0}$. It is well-known that any $\gamma\in \mcr{E}(\pi)$ has the unique decomposition $\gamma= \gamma^i + \gamma^p$ for $\gamma^i\in \mcr{E}_i(\pi)$ and $\gamma^p\in \mcr{E}_p(\pi)$; see, e.g., Getoor and Glover (1987). There is also a close connection between the classes $\mcr{E}^*_p(Q)$ and $\mcr{E}_p(\pi)$ involving immigration. In fact, to each $\gamma\in \mcr{E}_p(\pi)$ there corresponds a unique $\eta\in \mcr{K}(\pi)$ such that
 \beqlb\label{5.2}
\gamma(f) = \int_0^\infty \eta_s(f)\d s, \qquad f\in B(E).
 \eeqlb
By Proposition~8.7 in Li (2011), there is a unique $\kappa\in \mcr{K}(P)$ so that $\eta_t(f)= \pi_t(\kappa,f)$ for $t>0$ and $f\in B(E)^+$. We can define $N\in \mcr{E}^*(Q)$ by
 \beqlb\label{5.3}
L_N(f)
 =
\exp\bigg\{-\int_0^\infty V_s(\kappa,f)\d s\bigg\}, \qquad f\in B(E)^+.
 \eeqlb
It is not hard to see that
 \beqlb\label{5.4}
\gamma(f) =\int_{M(E)^\circ}\nu(f)N(\d\nu), \qquad f\in B(E).
 \eeqlb
Those relations establish a connection between the classes $\mcr{E}(P)$ and $\mcr{E}^*(Q)$. The connections between C-excessive distributions and excessive measures for the transition semigroup $(Q_t^\circ)_{t\ge 0}$ were discussed in Li (2003, 2011).

Let $(g_t)_{t\ge 0}$ be the composition semigroup of probability generating functions of a continuous-time branching process. A probability generating function $f$ is called \textit{self-decomposable} relative to $(g_t)_{t\ge 0}$ by Van Harn et al.\ (1982) if for each $t\ge 0$ there is another probability generating function $f_t$ so that
 \beqlb\label{5.5}
f(z) = (f\circ g_t)(z) f_t(z), \qquad |z|\le 1.
 \eeqlb
This generalizes the classical concept of self-decomposability; see, e.g., Lo\`eve (1977) and Sato (1999). A general representation of self-decomposable probability generating functions for a critical or subcritical branching process was given in Van Harn et al.\ (1982); see also the earlier work of Steutel and Van Harn (1979). In view of \eqref{5.1} and \eqref{5.5}, we may regard \eqref{4.21} as a counterpart of the representation (6.1b) of Van Harn et al.\ (1982) in the setting of Dawson--Watanabe superprocesses.

\section{Examples}

 \setcounter{equation}{0}

\bgexample\label{e6.1} Let $\xi$ be a Borel right process in $E$ with transition semigroup $(P_t)_{t\ge 0}$ and $\phi_1$ a local branching mechanism given by \eqref{3.14} with $\gamma(\cdot,1)\equiv 0$. Let $m(x,\d u)$ be the image of $H(x,\d\nu)$ under the mapping $\nu\mapsto \nu(1)$. Then $(u\land u^2)m(x,\d u)$ be a bounded kernel from $E$ to $(0,\infty)$ and the local projection $\phi_1$ has the representation
 \beqlb\label{6.1}
\phi_1(x,z)= b(x)z + c(x)z^2 + \int_{(0,\infty)} (\e^{-zu} - 1 + zu) m(x,\d u), \quad x\in E, z\ge 0.
 \eeqlb
In this case, the cumulant semigroup of the $(\xi,\phi_1)$-superprocess is defined by
 \beqnn
V_tf(x) = P_tf(x) - \int_0^t\d s\int_E \phi_1(y,V_sf(y)) P_{t-s}(x,\d y), \quad
x\in E, t\ge 0.
 \eeqnn
For any $\eta\in M(E)$, we can define the transition semigroup $(Q^\eta_t)_{t\ge 0}$ of an immigration superprocess by
 \beqlb\label{6.2}
\int_{M(E)} \e^{-\nu(f)} Q^\eta_t(\mu,\d\nu)
 =
\exp\bigg\{-\mu(V_tf)-\int_0^t\eta(V_sf)\d s\bigg\}, \qquad f\in B(E)^+.
 \eeqlb
Let $b_*= \inf_{x\in E}b(x)$ and $c_*= \inf_{x\in E}c(x)$. From our general results, we derive immediately the following properties of the process:
 \bitemize

\itm We have $L_{\var}(Q^\eta_tF)\le e^{-b_*t}L_{\var}(F)$ for $t\ge 0$ and Borel function $F$ on $M(E)$.

\itm If $c_*> 0$, then $L_{\var}(Q^\eta_tF)\le 2\|\bar{V}_t\|\|F\|$ for $t> 0$ and Borel function $F$ on $M(E)$.

\itm If $b_*> 0$, then $(Q^\eta_t)_{t\ge 0}$ has a unique stationary distribution $N_\infty$ and there is a constant $C\ge 0$ so that
 \beqnn
W_1(Q^\eta_t(\mu,\cdot),N_\infty)\le C(1+\mu(1))e^{-b_*t}, \quad t\ge 0, \mu\in M(E).
 \eeqnn

\itm If $c_*> 0$ and $b_*> 0$, there is a constant $C\ge 0$ so that
 \beqnn
\|Q^\eta_t(\mu,\cdot)-N_\infty\|_{\var}\le C(1+\mu(1))e^{-b_*t}, \quad t\ge 0, \mu\in M(E).
 \eeqnn

 \eitemize
The above results generalize those of Stannat (2003b, Theorems~1.7, 2.5 and 3.1), who assumed $E$ is a compact metric space, $\xi$ is a Feller process and $(x,z)\mapsto \phi_1(x,z)$ is jointly continuous on $E\times [0,\infty)$; see also Stannat (2003a). \edexample

\bgexample\label{e6.2} Let $0<\alpha<1$ and $a\in B(E)^+$. A typical special form of \eqref{6.1} is $\phi_1(x,z)= b(x)z + c(x)z^2 + a(x)z^{1+\alpha}$. This branching mechanism was excluded by the results in Stannat (2003a, 2003b). Let $a_*= \inf_{x\in E} a(x)$ and define $b_*$ and $c_*$ similarly from $b\in B(E)$ and $c\in B(E)^+$. Then Condition~\ref{t3.3} holds with $\phi_*(z)= b_*z + c_*z^2 + a_*z^{1+\alpha}$, which satisfies Grey's condition \eqref{3.8} if and only if $c_*+a_*> 0$. \edexample

\bgexample\label{e6.3} Consider a Borel right underlying process $\xi$ in $E$ with transition semigroup $(P_t)_{t\ge 0}$ and a branching mechanism $\phi$ given by \eqref{3.1} or \eqref{3.2}. Let $(V_t)_{t\ge 0}$ be the cumulant semigroup defined by \eqref{3.3}. Suppose that $\eta\in M(E)$ and $\nu(1)H(\d\nu)$ is a finite measure on $M(E)^\circ$. For $f\in B(E)^+$ write
 \beqnn
I(\eta,H,f) = \eta(f) + \int_{M(E)^\circ} \big(1-\e^{-\nu(f)}\big) H(\d \nu).
 \eeqnn
By Theorem~\ref{t4.13}, we can define the transition semigroup $(Q^N_t)_{t\ge 0}$ of an immigration superprocess by
 \beqlb\label{6.3}
\int_{M(E)} \e^{-\nu(f)} Q^N_t(\mu,\d\nu)
 =
\exp\bigg\{-\mu(V_tf)-\int_0^tI(\eta,H,V_sf)\d s\bigg\}.
 \eeqlb
This is a special case of \eqref{4.19} and a generalization of \eqref{6.2}. The immigration structure is determined by a closable infinitely divisible probability entrance law via \eqref{4.3}. The ergodicity and exponential ergodicity of this semigroup in the Wasserstein distance and another distance defined by Laplace functionals were studied in Friesen (2019+), where the state space was enlarged as described in Remark~\ref{t4.20} to include some tempered measures. \edexample

\bgexample\label{e6.4}
Suppose that $E$ is a bounded domain in $\mbb{R}^d$ with twice continuously differentiable boundary $\partial E$. Let $\xi$ be an absorbing-barrier Brownian motion in $E$ with transition semigroup $(P_t)_{t\ge 0}$. Let $(V_t)_{t\ge 0}$ be the cumulant semigroup defined by \eqref{3.3} for a branching mechanism $\phi$ given by \eqref{3.1} or \eqref{3.2}. It is well-known that $P_t(x,\d y)$ has a symmetric density $p_t(x,y) = p_t(y,x)$ for $t>0$, which is the fundamental solution of the heat equation on $E$ with Dirichlet boundary condition. Moreover, the density $p_t(x,y)$ is continuously differentiable in $x$ and $y$ to the boundary $\partial E$; see, e.g., Friedman (1964, p.83). Then for any $t>0$ and $f\in B(E)$ the function $P_tf$ is smooth on $E$ and we can extend it trivially to $\partial E$ by continuity. Let $h$ be the bounded strictly positive excessive function for $(P_t)_{t\ge 0}$ defined by
 \beqlb\label{6.4}
h(x) = \int_0^1P_s1(x)\d s, \qquad x\in E.
 \eeqlb
Then $h(x)\to 0$ as $x\to z\in \partial E$. Let $M_h(E)$ denote the set of $\sigma$-finite measures $\mu$ on $E$ such that $\mu(h)< \infty$. We also use $\partial$ to denote the operator of inward normal differentiation at the boundary $\partial E$. By Theorem~8.26 of Li (2011), any entrance law $\kappa\in \mcr{K}(P)$ has the representation
 \beqnn
\kappa_t(f) = \eta(P_tf) + \gamma(\partial P_tf), \quad t>0, f\in B(E),
 \eeqnn
where $\eta\in M_h(E)$ and $\gamma\in M(\partial E)$. It is not hard to show that
 \beqnn
V_t(\kappa,f) = \eta(V_tf) + \gamma(\partial V_tf), \quad t>0, f\in B(E)^+.
 \eeqnn
Suppose that $F(\d\nu,\d\zeta)$ is a $\sigma$-finite measure on $M(E)\times M(\partial E)$ satisfying
 \beqnn
\int_{M_h(E)\times M(\partial E)} \big[\nu(h) + \zeta(\partial h)\big] F(\d\nu,\d\zeta)< \infty.
 \eeqnn
By Theorem~\ref{t4.13}, the transition semigroup $(Q^N_t)_{t\ge 0}$ of an immigration process associated with the super absorbing-barrier Brownian motion has the representation
 \beqlb\label{6.5}
\int_{M(E)} \e^{-\nu(f)} Q^N_t(\mu,\d\nu)
 =
\exp\bigg\{-\mu(V_tf)-\int_0^t I_s(\kappa,F,f)\d s\bigg\},
 \eeqlb
where
 \beqnn
I_s(\kappa,F,f) = \eta(V_sf) + \gamma(\partial V_sf) + \int_{M_h(E)}\int_{M(\partial E)} \big(1-\e^{-\nu(V_sf) + \zeta(\partial V_sf)}\big) F(\d\nu,\d\zeta).
 \eeqnn
This is not a special case of the transition semigroup defined by \eqref{6.3} unless $\eta\in M(E)$, $\gamma(\partial E)=0$ and $F$ is supported by $M(E)\times \{0\}$.
\edexample

\bgexample\label{e6.5} Let $\xi$ be the standard absorbing-barrier Brownian motion in $E_0= (0,\infty)$ and let $(P_t)_{t\ge 0}$ denote its transition semigroup. For $t>0$ the kernel $P_t(x,\d y)$ has density
 \beqnn
p_t(x,y) = g_t(x-y) - g_t(x+y), \quad x,y>0,
 \eeqnn
where
 \beqnn
g_t(z) = \frac{1}{\sqrt{2\pi t}} \exp\{-z^2/2t\}, \quad t>0, z\in\mbb{R}.
 \eeqnn
Let $h$ be the bounded strictly positive excessive function for $(P_t)_{t\ge 0}$ defined as in Example~\ref{e6.4}. Let $M_h(E_0)$ denote the set of $\sigma$-finite measures $\mu$ on $E_0$ such that $\mu(h)< \infty$. Let $\partial$ denote the operator of upward normal differentiation. By Theorem~8.28 of Li (2011), an entrance law $\kappa\in \mcr{K}(P)$ has the representation
 \beqnn
\kappa_t(f) = \eta(P_tf) + a\partial P_tf(0), \quad t>0, f\in B(E_0),
 \eeqnn
where $\eta\in M_h(E_0)$ and $a\in [0,\infty)$. Suppose that $F(\d\nu,\d z)$ is a $\sigma$-finite measure on $M_h(E_0)\times [0,\infty)$ satisfying
 \beqnn
\int_{M_h(E_0)\times [0,\infty)} [\nu(h) + z] F(\d\nu,\d z)< \infty.
 \eeqnn
By Theorem~\ref{t4.13}, an immigration process associated with the super absorbing-barrier Brownian motion has transition semigroup $(Q^N_t)_{t\ge 0}$ defined by
 \beqlb\label{6.6}
\int_{M(E_0)} \e^{-\nu(f)} Q^N_t(\mu,\d\nu)
 =
\exp\bigg\{-\mu(V_tf)-\int_0^t I_s(\kappa,F,f)\d s\bigg\},
 \eeqlb
where
 \beqnn
I_s(\kappa,F,f) = \eta(V_sf) + a\partial V_sf(0) + \int_{M_h(E_0)}\int_{[0,\infty)} \big(1-\e^{-\nu(V_sf) + z\partial V_sf(0)}\big) F(\d\nu,\d z).
 \eeqnn
This is not a special case of the transition semigroup defined by \eqref{6.3} unless $\eta\in M(E)$, $a=0$ and $F$ is supported by $M(E_0)\times \{0\}$. The reader may refer to Section~9.4 of Li (2011) for explanations of the immigration from the origin involved in the semigroup.  \edexample

\bgexample\label{e6.6} Let $\langle\cdot,\cdot\rangle$ denote the Euclidean inner product of $\mbb{R}^d$. For each $i=1,\ldots,d$ let $\phi_i$ be a function on $\mbb{R}_+^d$ given by
 \beqlb\label{6.7}
\phi_i(\lambda) = b_i\lambda_i + c_i\lambda_i^2 - \langle\eta_i,\lambda\rangle +
\int_{\mbb{R}_+^d\setminus \{0\}} \big(\e^{-\langle\lambda,u\rangle} - 1 +
\lambda_iu_i\big) H_i(\d u),
 \eeqlb
where $c_i\ge 0$ and $b_i$ are constants, $\eta_i= (\eta_{i1},\cdots,\eta_{id})\in \mbb{R}_+^d$ is a vector with $\eta_{ii}=0$ for $i=1,\ldots,d$, and $H_i(\d u)= H_i(\d u_1,\cdots,\d u_d)$ is a $\sigma$-finite measure on $\mbb{R}_+^d\setminus \{0\}$ so that
 \beqnn
\int_{\mbb{R}_+^d\setminus \{0\}} \big(\langle u,1\rangle\land \langle u,1\rangle^2 + \langle u,1\rangle - u_i\big) H_i(\d u)< \infty.
 \eeqnn
For any $\lambda\in \mbb{R}_+^d$ there is a unique locally bounded vector-valued
solution $t\mapsto v(t,\lambda)\in \mbb{R}_+^d$ to the evolution equation system
 \beqlb\label{6.8}
\frac{\d v_i}{\d t}(t,\lambda) = - \phi_i(v(t,\lambda)), \quad
v_i(0,\lambda) = \lambda_i, \qquad i=1,\ldots,d.
 \eeqlb
We can define a transition semigroup $(Q_t)_{t\ge 0}$ on $\mbb{R}_+^d$ by
 \beqlb\label{6.9}
\int_{\mbb{R}_+^d}\e^{-\langle\lambda,y\rangle}Q_t(x,\d y)
 =
\exp\{-\langle x,v(t,\lambda)\rangle\}, \qquad \lambda, x\in \mbb{R}_+^d.
 \eeqlb
By the result of Rhyzhov and Skorokhod (1970), up to a moment assumption, this gives the most general form of a stochastically continuous transition semigroup on $\mbb{R}_+^d$ satisfying the branching property. A Markov process in $\mbb{R}_+^d$ with transition semigroup $(Q_t)_{t\ge 0}$ given by \eqref{6.9} is called a \textit{multi-type CB-process}. Let $\psi$ be a function on $\mbb{R}_+^d$ with the representation
 \beqlb\label{6.10}
\psi(\lambda)= \langle\beta,\lambda\rangle + \int_{\mbb{R}_+^d\setminus \{0\}} \big(1 - \e^{-\langle\lambda,u\rangle}\big) \nu(\d u),
 \eeqlb
where $\beta\in \mbb{R}_+^d$ is a vector and $\nu(\d u)= \nu(\d u_1,\cdots,\d u_d)$ is a $\sigma$-finite measure on $\mbb{R}_+^d\setminus \{0\}$ so that
 \beqnn
\int_{\mbb{R}_+^d\setminus \{0\}} \langle u,1\rangle\nu(\d u)< \infty.
 \eeqnn
We can define another transition semigroup $(Q^N_t)_{t\ge 0}$ on $\mbb{R}_+^d$ by
 \beqlb\label{6.11}
\int_{\mbb{R}_+^d}\e^{-\langle\lambda,y\rangle}Q^N_t(x,\d y)
 =
\exp\bigg\{-\langle x,v(t,\lambda)\rangle - \int_0^t\psi(v(s,\lambda))\d s\bigg\}, \quad \lambda, x\in \mbb{R}_+^d.
 \eeqlb
This is the finite-dimensional version of \eqref{6.3}. A Markov process in $\mbb{R}_+^d$ with transition semigroup $(Q^N_t)_{t\ge 0}$ given by \eqref{6.11} is called a \textit{multi-type CBI-process}. The process has been used widely in mathematical finance as models for interest rates or asset prices; see, e.g., Duffie et al.\ (2003). Let
 \beqnn
\gamma_{ij}= \eta_{ij} + \int_{\mbb{R}_+^d\setminus \{0\}} 1_{\{i\neq j\}}u_j H_i(\d u), \qquad i,j= 1,\cdots,d.
 \eeqnn
Now suppose that there is a branching mechanism $\phi_*$ in the form \eqref{3.5} satisfying Grey's condition \eqref{3.8} so that
 \beqnn
(b_i-\langle\gamma_i,1\rangle)z + c_iz^2 + \int_{\mbb{R}_+^d\setminus \{0\}} \big(\e^{-zu_i} - 1 + zu_i\big) H_i(\d u)\ge \phi_*(z), \qquad z\ge 0.
 \eeqnn
By Corollaries~\ref{t3.6} and~\ref{t4.4}, the semigroup $(Q^N_t)_{t\ge 0}$ has the strong Feller property. This generalizes the second assertion of Corollary~3.2 in Stannat (2003b), who considered the situation where $\nu=0$ and $H_i(\d u)$ is carried by the $i$th half-axis $\{u= (u_1,\ldots,u_d)\in \mbb{R}_+^d: u_i>0$ and $u_j=0$ if $j\neq i\}$ for $i=1,\ldots,d$. By Corollary~\ref{t4.18}, if $\beta_*:= \min_{1\le i\le d} (b_i - \langle\gamma_i, 1\rangle)> 0$, the multi-type CBI-process is exponentially ergodic in the total variation distance.
 \edexample

\bgexample\label{e6.7} Suppose that $\phi_*$ is the local branching mechanism given by \eqref{3.5} and $(v_t^*)_{t\ge 0}$ is defined by \eqref{3.7}. For $t\ge 0$ let $g_t(x)= (x-t)\vee 0$ if $x\ge 0$ and $= (x+t)\land 0$ if $x<0$. Let $\xi= \{\xi_t: t\ge 0\}$ be the Markov process in $\mbb{R}$ satisfying $\xi_t= g_t(\xi_0)$ for $t\ge 0$. Suppose that $X = (W, \mcr{G}, \mcr{G}_t, X_t, \mbf{Q}_\mu)$ is a right realization of the $(\xi,\phi_*)$-superprocess. By a modification of the proof of Proposition~5.20 in Li (2011), one can show that $\mbf{Q}_{a\delta_x}\{X_t = X_t(1) \delta_{g_t(x)}$ for $t\ge 0\} = 1$ for $x\in \mbb{R}$ and $a\ge 0$. The cumulant semigroup $(V_t)_{t\ge 0}$ of this $(\xi,\phi_*)$-superprocess is given by $V_tf(x)= v^*_{t\land |x|}(g_t(x))$ for $x\in \mbb{R}$ and $t\ge 0$. Suppose that $\mu_1\in M(\mbb{R})$ has bounded support $\supp(\mu_1)\subset [0,c]$, where $c\ge 0$. Let $\mu_2$ be the image of $\mu_1$ induced by the mapping $x\mapsto -x$. For $t\ge c$ and $f\in B(\mbb{R})^+$ we have
 \beqnn
\int_{M(\mbb{R})}\e^{-\nu(f)}Q_t(\mu_i,\d\nu)
 =
\exp\bigg\{-\int_{[0,c]}v_x^*(f(0))\mu_1(\d x)\bigg\}, \quad i=1,2,
 \eeqnn
and hence $Q_t(\mu_1,\cdot)= Q_t(\mu_2,\cdot)$. It follows that
 \beqnn
W_1(Q_t(\mu_1,\cdot), Q_t(\mu_2,\cdot)) = \|Q_t(\mu_1,\cdot)- Q_t(\mu_2,\cdot)\|_{\var} = 0.
 \eeqnn
Then the lower bounds given in Theorems~\ref{t2.2} and~\ref{t2.8} for this $(\xi,\phi_*)$-superprocess can be reached. \edexample

\bgexample\label{e6.8} Suppose that $\phi_*$ and $(v_t^*)_{t\ge 0}$ are given as in the last example. Let $\xi$ be the Markov process in $[0,1]$ defined by $\xi_t= (\xi_0-t)\vee \lfloor\xi_0\rfloor$ for $t\ge 0$, where ``$\lfloor\cdot\rfloor$'' denote the integer part. Let $\phi$ be the local branching mechanism on $[0,1]$ defined by $\phi(x,z) = 1_{(0,1]}(x)\phi_*(z)$ for $x\in [0,1]$ and $z\ge 0$. Suppose that $X = (W, \mcr{G}, \mcr{G}_t, X_t, \mbf{Q}_\mu)$ is a right realization of the $(\xi,\phi)$-superprocess. Then $\mbf{Q}_{a\delta_x}\{X_t = X_t(1)\delta_{(x-t)\vee \lfloor x\rfloor}$ for $t\ge 0\} = 1$ for $x\in [0,1]$ and $a\ge 0$. The cumulant semigroup $(V_t)_{t\ge 0}$ of this $(\xi,\phi)$-superprocess is given by $V_tf(x)= v^*_{t\land x}(f((x-t)\vee 0))$ for $x\in [0,1)$ and $V_tf(x)= v_t^*(f(1))$ for $x=1$. Clearly, for any $\mu\in M([0,1))\subset M([0,1])$, the limit $Q_\infty(\mu,\cdot):= \lim_{t\to \infty}Q_t(\mu,\cdot)$ exists and, for $f\in B([0,1])^+$,
 \beqnn
\int_{M([0,1])}\e^{-\nu(f)}Q_\infty(\mu,\d\nu)
 =
\exp\bigg\{-\int_{[0,1)}v_z^*(f(0))\mu(\d z)\bigg\}.
 \eeqnn
It is easy to see that $Q_\infty(\mu,\cdot)\in \mcr{E}^*_i(Q)$ is carried by $M(\{0\})\subset M([0,1])$. If $b_*= \phi_*^\prime(0)> 0$, for each $\beta>0$ we can define $N^\beta\in \mcr{E}_p^*(Q)$ by, for $f\in B([0,1])^+$,
 \beqnn
\int_{M([0,1])}\e^{-\nu(f)}N^\beta(\d\nu)
 =
\exp\bigg\{-\beta\int_0^\infty v_s^*(f(1))\d s\bigg\}.
 \eeqnn
Then both $\mcr{E}_i^*(Q)$ and $\mcr{E}_p^*(Q)$ contain non-trivial elements.
\edexample

\bigskip

\textbf{Acknowledgments} ~ I would like to thank Professors Yonghua Mao and Yuhui Zhang for helpful discussions on Wasserstein distances. I am grateful to the Laboratory of Mathematics and Complex Systems (Ministry of Education) for providing the research facilities to carry out the project.


 \end{document}